\documentclass[11pt]{article}
\usepackage{amstex}
\usepackage{amssymb}
\font\tenmsam=msam10
\def\square{\hbox{\tenmsam\char'003}}

\def\hfl#1#2{\smash{\mathop{\hbox to 12mm{\rightarrowfill}}
\limits^{\scriptstyle#1}_{\scriptstyle#2}}}

\newtheorem{theo}{Theorem}[subsection]
\newtheorem{prop}[theo]{Proposition}
\newtheorem{coro}[theo]{Corollary}
\newtheorem{lemm}[theo]{Lemma}
\newtheorem{defi}[theo]{Definition}
\newtheorem{rema}[theo]{Remark}
\newtheorem{ex}[theo]{Example}

\title{\textsc{galois extension for a compact quantum group}}
\date{} 
\author{Julien Bichon}

\makeatletter
\renewcommand{\@makefnmark}{}
\makeatother

\begin{document}
\maketitle

\begin{abstract}
The aim of this paper is to introduce the quantum analogues of
torsors for a compact quantum group
and to investigate their relations with representation theory.
Let $A$ be a Hopf algebra over a field $k$. A theorem of Ulbrich asserts that
there is an equivalence of categories between neutral fibre
functors on the category of finite-dimensional $A$-comodules and 
$A$-Galois extensions of $k$. We give the compact quantum group
version of this result. Let $A$ be the Hopf $*$-algebra
 of representative functions on a
compact quantum group. We show that there is an equivalence of categories
between $*$-fibre functors on unitary A-comodules and 
$A$-$*$-Galois extensions with a positive
Haar measure. Such an  $A$-$*$-Galois extension has a $C^*$-norm, which
furthermore can be taken as the upperbound of $C^*$-semi-norms. We then
introduce the notion of Galois extension for a 
compact quantum group, for which measure theory can be deduced from topology.
We  construct universal Galois extensions, 
which enable us to find a nontrivial Galois extension 
for the unitary quantum group  $U_q(2)$.
\end{abstract}

\footnote{\noindent Département de Mathématiques, Université
Montpellier II. \\
 Place Eugène Bataillon, 34095 Montpellier Cedex 5.
E-mail address: bichon\char64 math.univ-montp2.fr}

\section{Introduction}

  A torsor for a group is a space on which the group 
acts freely and transitively. There is an obvious group 
structure on a torsor. This is no more true 
for the non-commutative analogue 
of torsor: Galois extensions for Hopf algebras.

In this paper, we introduce and investigate the notion 
of Galois extension for a compact quantum group. We mostly 
will be concerned with the tannakian use of 
torsors: classification of fibre functors on the category 
of representations of a compact quantum group. 
Moreover, Galois extensions provide new and interesting 
examples of quantum spaces.

\medskip

Let $k$ be a commutative field and let $G$ be a $k$-affine group scheme. 
We know from \cite{[S]} or \cite{[DM]} that there is an equivalence of 
categories  between symmetric fibre functors on the category of 
representations of $G$ and $G$-torsors of $k$.

The notion of Galois extension for an arbitrary Hopf $k$-algebra $A$ is 
the non-commutative version of torsor: an $A$-Galois extension of $k$ 
is a left $A$-comodule algebra for which a certain 
map $Z \otimes Z \longrightarrow A \otimes Z$ is an isomorphism. 
Ulbrich (\cite{[U1], [U2]}) then constructs an equivalence of 
categories between fibre functors on $A$-comodules 
and $A$-Galois extensions of $k$.

\medskip

We give here a compact quantum group version of Ulbrich's theorem. 
Two Hopf algebras are said to be Hopf co-Morita equivalent if there 
is a monoidal equivalence between their categories of finite-dimensional 
comodules. Let $A$ be a Hopf algebra. We naturally would like to find all 
Hopf algebras Hopf co-Morita equivalent to $A$. We know from tannakian 
duality that it is enough to describe the category of fibre functors 
on finite-dimensional $A$-comodules for this purpose. Ulbrich's theorem, 
which gives a convenient technique to construct fibre functors, 
and ensures that all fibre functors arise in this way, is thus very 
important in the Hopf co-Morita theory setting. Our version of his 
theorem will have the same role in the Hopf co-Morita theory of 
compact quantum group (see 4.5 for the precise definition). 
A good understanding of Hopf co-Morita theory
is certainly a necessary step towards an abstract duality for
compact quantum groups in the spirit of Doplicher and Roberts (\cite{[DR]}). 

\medskip

In this paper, we use the following  definition (5.3) 
(equivalent to the classical ones \cite{[W3], [W1], [DK], [MV]}):
 a compact quantum group is a pair  $(A, ||.||)$ where $A$
is a Hopf $*$-algebra and $|| . ||$ is a $C^*$-norm on
$A$ such that the comultiplication is continuous. The Hopf $*$-algebra 
$A$ is called the algebra of representative 
functions on the compact quantum group. 

A Hopf $*$-algebra is called unitarizable if every $A$-comodule is isomorphic 
to a unitary one. One of the main results in quantum group theory 
(mainly due to Woronowicz, but also to Van Daele  in full generality) is 
that the Hopf $*$-algebra of representative functions on a compact 
quantum group is unitarizable (\cite{[V], [W1], [W3]}). 
Conversely, any unitarizable 
Hopf $*$-algebra gives rise to a
compact quantum group (see \cite{[DK],[W2]}, a unitarizable 
Hopf $*$-algebra  is a CQG algebra in the sense of \cite{[DK]}).

\smallskip
     
Let $A$ be a unitarizable Hopf $*$-algebra. Denote by $\mathcal U_f(A)$ 
the monoidal $*$-category of unitary $A$-comodules. The category 
$\mathcal U_f(A)$
may be seen as the category of unitary representations of  
the quantum groups whose algebra of representative functions is $A$.
A $*$-fibre functor on  
 $\mathcal U_f(A)$ is a monoidal $*$-functor  
$\eta : \mathcal U_f(A) \longrightarrow \mathcal H$
(the category of finite-dimensional Hilbert spaces) 
whose monoidal functor constraints 
$\eta (V) \otimes \eta(W)\ \widetilde {\longrightarrow}\
\eta(V \otimes W)$ and  
$\mathbb C \ \widetilde{\longrightarrow}\ \eta(\mathbb C)$
are unitaries. A $*$-fibre functor is a special fibre functor.

We show that there is an equivalence of categories between $*$-fibre functors
on  $\mathcal U_f(A)$ and $A$-$*$-Galois extension with a positive Haar measure,
by using an adaptation of Ulbrich's functor \cite{[U2]}. Our result, although
related to representations of compact quantum groups, is stated in a 
purely ($*$)-algebraic way. This is not a surprise. 
Indeed, the reconstructed object in the Woronowicz-Tannaka-Krein theorem 
for compact quantum groups (\cite{[W2]})
is the algebra of representative functions on the compact quantum group 
with its maximal $C^*$-norm, not the a priori given $C^*$-norm.

Next, we construct a $C^*$-norm on an  $A$-$*$-Galois extension with 
a positive Haar measure, and show that the upperbound of 
$C^*$-semi-norms exists.
The notion of Galois extension for a compact quantum group is then 
introduced, without postponing the existence of a Haar measure. 
We construct the positive Haar measure on these Galois extensions.

Using our results, we find another proof of 
a theorem from Doplicher-Roberts:
a symmetric $*$-fibre functor on the category of 
unitary representations of a compact group is unique up to unitary 
isomorphism (\cite{[DR]}).

The universal $C^*$-algebras $O_{n,p}$ generated by an 
unitary $(n,p)$-matrix arise naturally in the context of 
Galois extensions (in the case $p=1$, one gets the Cuntz algebras 
\cite{[C]}). They are closely related to 
the property for a $*$-fibre functor to preserve dimensions of the underlying 
vector spaces.

We also define universal Galois extensions (in the spirit of Van Daele and Wang \cite{[VW]}), useful to show the existence of Haar measures. They allow us to find a non trivial Galois extension for the quantum group $U_q(2)$
(unitary quantum group), and hence a non trivial $*$-fibre functor on its 
category of unitary representations.

\medskip

Our work is organized as follows. In section 2, we review the material on 
Galois extensions and Ulbrich's theorem. The inverse functor is described in a 
different way from \cite{[U2]}. Our inverse functor is inspired from the 
tannakian theorem \cite{[U3]}. We find a non-trivial Galois extension for 
$GL_q(2)$, which will be used to construct a Galois extension
 for $U_q(2)$ in section 6.
In section 3, we review some  facts on  unitary comodules, $*$-categories 
and $*$-fibre functors. 
Section 4 is devoted to the equivalence of categories between Galois extensions and $*$-fibre functors.
In section 5, We construct $C^*$-norms on our Galois extensions and we define the notion of Galois extension for a compact quantum group. We also give our 
proof of Doplicher-Roberts' unicity theorem.
In section 6, we define universal Galois extensions.

\subsection*{Notations}
All algebras will have units.

Let  $A = (A,m,u,\Delta,\varepsilon,S)$ be a Hopf $k$-algebra.
 The multiplication will be denoted by  $m$,   $u : k
\rightarrow A$ is the unit of $A$, while $\Delta$, $\varepsilon$ and $S$
are respectively the  comultiplication, the counit and the antipode of $A$.

Let us denote by ${\rm Co}_f(A)$ the category of finite-dimensional 
(right) $A$-comodules. When $V$ is a (right) $A$-comodule, we denote by 
$\alpha_V : V \longrightarrow V \otimes A$
the coaction. The category of finite dimensional $k$-vector spaces will
be denoted by ${\rm Vect}_f(k)$.
Galois extensions will be left $A$-comodule algebras.
By an $A$-Galois extension we will mean a left $A$-Galois extension of $k$
 as in \cite{[Sc],[U1],[U2]}.

We refer the reader to \cite{[B], [JS]} for the definitions of 
monoidal categories, monoidal functors  and for tannakian duality.

\medskip
\noindent
{\bf Aknowledgements}. The author wishes to thank
Guy Laffaille for his numerous advices and his patient reading
of the manuscript.  
\section{Galois extensions}

In the whole section $A$ will be a Hopf algebra over a commutative field
$k$. We will denote by  $\omega : {\rm Co}_f(A) \longrightarrow {\rm
Vect}_f(k)$ the forgetful functor on finite-dimensional $A$-comodules.

\subsection{Galois extensions and fibre functors}
Let $G$ be a group. A (left) $G$-space $X$ is said to be a $G$-torsor
if $G$ acts freely and transitively on $X$. The $G$-space $X$
is a $G$-torsor if an only if the following map is a bijection:
\begin{eqnarray*}
G \times X & \longrightarrow & X \times X \\
(g,x) & \longmapsto & (gx,x).
\end{eqnarray*}
This observation motivates the following definition
(the terminology comes from classical Galois theory, see \cite{[KT]}): 

\begin{defi}
A non-zero (left) $A$-comodule algebra $Z$ is said to be an 
{\rm $A$-Galois extension} if the linear map:
$$
(1_A \otimes m_Z) \circ (\alpha_Z \otimes 1_Z) : Z \otimes Z \longrightarrow A
\otimes Z
$$
is an isomorphism.

The {\rm category of $A$-galois extensions}, denoted {\rm Gal($A$)}, 
is the category whose objects are $A$-Galois extensions and whose 
morphisms are $A$-colinear algebras morphisms.
\end{defi}

\begin{ex}
{\rm The Hopf algebra $A$ is naturally an $A$-Galois extension,
with coaction defined by $\Delta : A \longrightarrow A \otimes A$.}
\end{ex}

\begin{defi}
A monoidal functor $\eta : {\rm Co}_f(A) \longrightarrow {\rm Vect}_f(k)$ is said to be a {\rm fibre functor} on  ${\rm Co}_f(A)$ if it is exact, $k$-linear and faithful.

The {\rm category of fibre functors} on ${\rm Co}_f(A)$, denoted {\rm Fib($A$)}, is the category whose objects are fibre functors on ${\rm Co}_f(A)$ and whose morphisms 
are morphisms of monoidal functors. The set of morphisms between fibre functors
$\eta$ and $\eta'$ will be denoted by ${\rm Hom}^\otimes (\eta,\eta')$. 
\end{defi}

\begin{rema}
{\rm A morphism between fibre functors is necessarily an isomorphism 
(see \cite{[B]} or \cite{[S]}).}
\end{rema}

\subsection{Ulbrich's theorem}

Let $Z$ be an $A$-Galois extension of $k$. Let us recall how a fibre 
functor is constructed from $Z$ (in \cite{[U1],[U2]}, Galois extensions 
are right comodule algebras and give rise to functors on left comodules).

Let us define $\eta_Z : {\rm Co}_f(A) \longrightarrow {\rm Vect}_f(k)$ by
$\eta_Z(V) = V \wedge Z$, where $V \wedge Z$ is the 
kernel of the double arrow:
$$
\alpha_V \otimes 1_Z,  \ 1_V \otimes \alpha_Z \ : \ V \otimes Z 
\rightrightarrows
V
\otimes A \otimes Z
$$
 and by $\eta_Z(f) = f \otimes 1_Z$
for $f \in {\rm Hom}_{{\rm Co}_f(A)}(V,W)$.
We know from \cite{[U1],[U2]} that $\eta_Z$ is a fibre functor. When
$Z =A$, then $\eta_A$ is isomorphic to the forgetful functor.
  For two comodule $V$ and $W$, the functorial isomorphism

 $$(V \wedge Z)
\otimes (W \wedge Z) \widetilde {\longrightarrow} (V \otimes W) \wedge Z$$
is induced by the linear map  $v \otimes z_1 \otimes w \otimes z_2
\longmapsto v \otimes w \otimes z_1z_2$.

If $\phi : Z \longrightarrow T$ is a morphism of $A$-Galois extensions, 
a morphism of fibre functors $\phi_* : \eta_Z \longrightarrow \eta_T$ 
is defined by letting $\phi_* = 1 \otimes \phi$.

We have just defined a functor ${\rm Gal}(A) \longrightarrow {\rm Fib}(A)$. The following result is proved in \cite{[U1],[U2]}:

\begin{theo}[Ulbrich]
 The functor  ${\rm Gal}(A) \longrightarrow {\rm Fib}(A)$, $Z
\longmapsto \eta_Z$, is an equivalence of categories.
\end{theo}

In subsections 2.3 and 2.4, we give a method to construct the quasi-inverse
functor. It is different from \cite{[U2]}: it is inspired by \cite{[U3]} and 
will be useful in the sequel.

\subsection{The Hom$^{\vee}$ construction}

Let $\eta : {\rm Co}_f(A) \longrightarrow {\rm Vect}_f(k)$ be a fibre functor,
and let $\omega$ be the forgetful functor.

Let $\mathcal N$ be the subspace of $\bigoplus_{V \in {\rm Co}_f(A)} {\rm
Hom}(\omega(V),\eta(V))$ generated by $\eta(f) \circ v - v \circ
\omega(f)$, with $f \in {\rm Hom}_{{\rm Co}_f(A)}(V,W)$ and $v \in {\rm
Hom}(\omega(W),\eta(V))$.

We define (\cite{[JS]}, sec. 3):
$$
{\rm Hom}^{\vee}(\eta,\omega) = \bigoplus_{V \in {\rm Co}_f(A)} {\rm
Hom}(\omega(V),\eta(V))/\mathcal N
.$$
If $v \in {\rm Hom}(\omega(V),\eta(V))$, we denote it by 
$[V,v]$ in  ${\rm Hom}^{\vee}(\eta,\omega)$.
Let  ${\rm Hom}^{\vee}(\omega,\omega) = {\rm End}^{\vee}(\omega) $. We know,
for example from \cite{[JS]}, that $A \cong {\rm End}^{\vee}(\omega)$ 
(tannakian reconstruction of $A$).  

Let us recall (\cite{[JS]}, sec. 4) that there is a linear map
$$
\delta : {\rm Hom}^{\vee}(\eta,\omega) \longrightarrow {\rm End}^{\vee}(\omega)
\otimes {\rm Hom}^{\vee}(\eta,\omega)
$$
which endows ${\rm Hom}^{\vee}(\eta,\omega)$ with a left
 ${\rm End}^{\vee}(\omega)$-comodule structure and hence
an $A$-comodule structure (since $A \widetilde
{\longrightarrow} {\rm End}^{\vee}(\omega) = {\rm Hom}^\vee(\omega, \omega)$). 
The map $\delta$ is defined in the following way:
  $\delta([V,\phi\otimes x])= \sum_i [V, \phi\otimes
v_i] \otimes [V, v_i^* \otimes x]$, with $V$ an
$A$-comodule with basis $v_1,\ldots,v_n$,  $\phi \in V^*$ and $x\in \eta(V)$.

Moreover  ${\rm Hom}^{\vee}(\eta, \omega)$ can be endowed with an algebra structure:
$$
[V,v].[W,w] = [V \otimes W, \widetilde \eta_{V,W} \circ (v \otimes w)]
$$
where $\widetilde \eta_{V,W} : \eta(V) \otimes \eta(W) \longrightarrow \eta(V
\otimes W)$  is the isomorphism given by the fibre functor $\eta$.
It is easy to check that ${\rm Hom}^{\vee}(\eta,\omega)$ 
is an ${\rm End}^{\vee}(\omega)$-comodule algebra.

Finally, let us see that ${\rm Hom}^{\vee}(\eta,\omega)$ is  ${\rm
End}^{\vee}(\omega)$-Galois extension. There is 
a linear map (\cite{[JS]}, sec. 4) $\delta' : {\rm End}^{\vee}(\omega) \longrightarrow {\rm
Hom}^{\vee}(\eta,\omega)
\otimes {\rm Hom}^{\vee}(\omega,\eta)$ and an ``antipode'' $S' : {\rm
Hom}^{\vee}(\omega,\eta)
\longrightarrow {\rm Hom}^{\vee}(\eta,\omega)$ induced by duality (cf. \cite{[U3]}
for the idea of construction ). The linear map:

\begin{eqnarray*}
\lefteqn{{\rm End}^{\vee}(\omega) \otimes {\rm Hom}^{\vee}(\eta,\omega)}\\
& & \stackrel{\delta'
\otimes 1}{\longrightarrow} {\rm Hom}^{\vee}(\eta,\omega) \otimes {\rm
Hom}^{\vee}(\omega,\eta) \otimes {\rm Hom}^{\vee}(\eta,\omega)  \\
& & \stackrel{1 \otimes S' \otimes 1}{\longrightarrow} {\rm
Hom}^{\vee}(\eta,\omega)
\otimes {\rm Hom}^{\vee}(\eta,\omega) \otimes {\rm Hom}^{\vee}(\eta,\omega) \\
& & \stackrel{1
\otimes m}{\longrightarrow} {\rm Hom}^{\vee}(\eta,\omega) \otimes {\rm
Hom}^{\vee}(\eta,\omega)
\end{eqnarray*}

\noindent
is inverse of the arrow (2.1.1) and hence ${\rm
Hom}^{\vee}(\eta,\omega)$ is a Galois extension.

If $\theta$ is another fibre functor and $u : \eta \longrightarrow \theta$
is a morphism of fibre functors, a morphism of $A$-Galois extensions 
$f :{\rm Hom}^{\vee}(\eta,\omega) \longrightarrow 
{\rm Hom}^{\vee}(\theta,\omega)$ is defined by 
$f([V,v])=[V,u_V \circ v]$.

We have just defined a functor 
${\rm Fib}(A) \longrightarrow {\rm Gal}(A)$, $\eta
\longmapsto {\rm Hom}^{\vee}(\eta,\omega)$.

\subsection{Sketch of proof}

Let us show that the functors 2.2 and 2.3 are quasi-inverse.
Let $\theta$ a fibre functor on  ${\rm Co}_f(A)$, and let $\eta_{{\rm
Hom}^{\vee}(\theta,\omega)}$ the fibre functor defined in 2.2.
There is (see \cite{[JS]}, sec. 4) a functorial morphism
  $\gamma : \theta
\longrightarrow \omega \otimes {\rm Hom}^{\vee}(\theta,\omega)$. 
It is easy to check that $\gamma$ gives rise to a morphism of fibre functors
$\theta
\longrightarrow
\eta_{{\rm Hom}^{\vee}(\theta,\omega)}$, which is an isomorphism 
(remark 2.1.4).

Let $Z$ be an $A$-Galois extension. The universal property of the
 ${\rm Hom}^{\vee}$ (\cite{[JS]}, sec. 4, prop. 3) gives rise to 
an extension morphism  ${\rm Hom}^\vee(\eta_Z,\omega) \longrightarrow Z$
which is an isomorphism.
(see \cite{[JS]}, sec. 6, minimal models for injectivity.
For surjectivity, one has to remark that every $z$ of $Z$ belongs to the 
range of the map $ev\otimes 1_Z : (V^*\otimes V)\wedge Z \rightarrow Z$, 
where $ev$ is the evaluation map, for some finite-dimensional $A$-comodule
$V$).

\subsection{The spectrum of a Galois extension}

Let $\eta : {\rm Co}_f(A) \longrightarrow {\rm Vect}_f(k)$ be a fibre functor. 
To any  $u \in {\rm Hom}(\eta,\omega)$ (natural transformation
from $\eta$ to $\omega$) is associated a linear form $f_u$ on
${\rm Hom}^\vee(\eta,\omega)$ in the following way  
: $f_u([V,v]) = {\rm Tr} (u_V \circ v)$. We get in this way
a linear isomorphism ${\rm Hom} (\eta,\omega)
\widetilde {\longrightarrow} {\rm Hom}^\vee(\eta,\omega)^*$  
(see \cite{[JS]}, sec. 3) which induces
another isomorphism:
$$
{\rm Hom}^\otimes (\eta,\omega) 
\widetilde {\longrightarrow} {\rm Hom}_{k-{\rm alg}} ({\rm
Hom}^\vee(\eta,\omega),k)
$$
where ${\rm Hom}^\otimes (\eta,\omega)$ is the set of morphisms of 
fibre functors $\eta \rightarrow \omega$.

\subsection{Hopf co-Morita equivalence}

\begin{defi}
Two Hopf $k$-algebras are said to be {\sl Hopf co-Morita equivalent} if their (right) finite-dimensional comodules categories are
monoidally equivalent.
\end{defi}

Let $\eta : {\rm Co}_f(A) \longrightarrow {\rm Vect}_f(k)$ be a fibre functor.
Let ${\rm End}^\vee(\eta)$ be the Hopf algebra created by tannakian duality
(constructed in the same way as in 2.3). 
We know from \cite{[JS]} or \cite{[B]} 
that $\eta$ factorizes through a monoidal equivalence
${\overline \eta} : {\rm Co}_f(A) 
{\widetilde {\longrightarrow}} {\rm Co}_f({\rm End}^\vee(\eta))$. 
This means that
 $A$ and  ${\rm End}^\vee(\eta)$ are Hopf co-Morita equivalent.
Moreover the Hopf algebras $A$ and ${\rm End}^\vee(\eta)$ are
isomorphic if and only if there is a monoidal equivalence 
$F : {\rm Co}_f(A) \rightarrow {\rm Co}_f(A)$ such that the fibre functors
$\eta \circ F$ and $\omega$ (forgetful functor on
${\rm Co}_f(A)$) are isomorphic. Conversely, any Hopf co-Morita equivalence
 for $A$ induces a fibre functor on ${\rm Co}_f(A)$.

A method to construct ${\rm End}^\vee(\eta)$
using the Galois extension associated to $\eta$ rather than $\eta$
itself is given in \cite{[Sc]}.

\subsection{A Galois extension for $GL_q(2)$}

The Hopf algebras $GL_{p,q}(2)$ and $GL_{p',q'}(2)$ are Hopf co-Morita
equivalent when $pq = p'q'$. In particular the Hopf algebras
 $GL_q(2)$ and $GL_{-q}(2)$ are Hopf co-Morita equivalent.
This fact is proved in \cite{[B]} with the use of quantum groupoids
(\cite{[M]}).
We find this result again by using Galois extensions.

\medskip
Let $k$ be a commutative field. Let us recall that $GL_q(2)$ (\cite{[K]}) 
is the quotient of the free algebra
$k\{x_{11},x_{12},x_{21},x_{22},t\}$ by the two-sided ideal generated
by the relations:
$$x_{12} x_{11} =  q x_{11} x_{12}, \quad
x_{21} x_{11} = q x_{11} x_{21}$$
$$x_{22} x_{12} = q x_{12} x_{22}, \quad
x_{22} x_{21} =  q x_{21} x_{22}$$
$$x_{12} x_{21} = x_{21} x_{12},\quad
x_{11} x_{22} - x_{22} x_{11} = (q^{-1} - q)x_{12}x_{21}$$
$$x_{11} t =  t x_{11},\quad
x_{22} t = t x_{22},\quad
x_{12} t =  t x_{12},\quad
x_{21} t =  t x_{21}$$
$$(x_{11} x_{22} - q^{-1}x_{12} x_{21}) t = 1$$

\begin{defi}
The algebra $GL_q(2,-2)$
is the universal algebra with generators
$z_{11}$, $z_{12}$, $z_{21}$, $z_{22}$, $t$ and relations:
$$z_{12} z_{11} = - q z_{11} z_{12}, \quad
z_{21} z_{11} = q z_{11} z_{21}$$
$$z_{22} z_{12} = q z_{12} z_{22}, \quad
z_{22} z_{21} = - q z_{21} z_{22}$$
$$z_{12} z_{21} = -z_{21} z_{12},\quad
z_{11} z_{22} + z_{22} z_{11} = (q-q^{-1})z_{12}z_{21}$$
$$z_{11} \tau = - \tau z_{11},\quad
z_{22} \tau = - \tau z_{22},\quad
z_{12} \tau = - \tau z_{12}, \quad
z_{21} \tau = - \tau z_{21}$$
$$(z_{11} z_{22} + q^{-1}z_{12} z_{21}) \tau = 1$$
\end{defi}

It is easy to check that  ${\rm Hom}_{k-{\rm alg}}(GL_q(2,-2),k) =
\emptyset$ and that $GL_q(2,-2)$ is a non-zero algebra.
 Direct computations lead to the following lemma:

\begin{lemm}
The algebra $GL_q(2,-2)$ is a left $GL_q(2)$-comodule algebra with coaction 
 $\alpha : GL_q(2,-2) \rightarrow
GL_q(2) \otimes GL_q(2,-2)$ defined by 
$\alpha(z_{ij}) = \sum_k x_{ik} \otimes
z_{kj}$ and $\alpha(\tau)= t\otimes \tau$.
\end{lemm}

We want to show that $GL_q(2,-2)$ is a $GL_q(2)$-Galois extension. 
For this purpose we introduce another algebra $GL_q(-2,2)$.

\begin{defi}
The algebra $GL_q(-2,2)$ is the universal algebra 
with generators $t_{11}$, $t_{12}$, $t_{21}$, $t_{22}$, $\xi$
and relations:
$$t_{12} t_{11} = qt_{11} t_{12}, \quad t_{21} t_{11} = -qt_{11}t_{21}$$
$$t_{22} t_{12} = - qt_{12} t_{22},\quad t_{22} t_{21} = qt_{21} t_{22}$$
$$t_{12}t_{21} = - t_{21}t_{12}, \quad
t_{11} t_{22} + t_{22} t_{11}= (q^{-1}-q)t_{12}t_{21}$$
$$t_{11} \xi = -\xi t_{11},\quad
t_{12} \xi = - \xi t_{12},\quad
t_{21} \xi = - \xi t_{21},\quad
t_{22} \xi= - \xi t_{22}$$
$$(t_{11} t_{22} - q^{-1}t_{12}t_{21}) \xi = 1$$
\end{defi}

Direct computations give the following two lemmas:

\begin{lemm}
 There is a morphism of algebras $\delta : GL_q(2,k)
\longrightarrow GL_q(2,-2) \otimes GL_q(-2,2)$ defined by
 $\delta(x_{ij}) = \sum_k
z_{ik} \otimes t_{kj}$ and $\delta(t)=\tau\otimes \xi$.
\end{lemm}

\begin{lemm}
The algebras $GL_q(2,-2)$ and $GL_q(-2,2)$ are anti-isomorphic
via the map $\phi : GL_q(-2,2)
\longrightarrow GL_q(2,-2)$ defined by $\phi(t_{11}) = z_{22} \tau$, $\phi(t_{12}) =
q\tau z_{12}$, $\phi(t_{21}) = q^{-1}z_{21} \tau$, $\phi(t_{22}) = \tau z_{11}$
and $\phi(\xi) = z_{11} z_{22} + q^{-1}z_{12} z_{21}$.
\end{lemm}

Let us define a linear map 
$GL_q(2) \otimes GL_q(2,-2) \longrightarrow GL_q(2,-2) \otimes
GL_q(2,-2)$ to be the composition of the maps:

\begin{eqnarray*}
& GL_q(2) \otimes GL_q(2,-2) & \stackrel{\delta \otimes 1}{\longrightarrow}
GL_q(2,-2)
\otimes GL_q(-2,2)\otimes GL_q(2,-2) \\
& & \stackrel{1 \otimes \phi \otimes 1}{\longrightarrow} GL_q(2,-2)
\otimes GL_q(2,-2) \otimes GL_q(2,-2) \\
& & \stackrel{1 \otimes m}{\longrightarrow} GL_q(2,-2)
\otimes GL_q(2,-2)
\end{eqnarray*}

It is easily seen that this map in inverse from the map 2.1.1. 
Indeed the  matrix \break $z = \left(\begin{array}{cc} z_{11} & z_{12} \\
z_{21} & z_{22} \end{array} \right)$ is invertible and its inverse 
is the matrix: 
$$ \left(\begin{array}{cc} \phi(t_{11}) & \phi(t_{12}) \\
\phi(t_{21}) & \phi(t_{22}) \end{array} \right) = \left( \begin{array}{cc}
z_{22} \tau & q\tau z_{12} \\ q^{-1}z_{21} \tau & \tau z_{11} \end{array} \right
).$$

Hence the algebra $GL_q(2,-2)$ is $GL_q(2)$-Galois extension and is 
non-trivial since ${\rm
Hom}_{k-{\rm alg}}(GL_q(2,-2),k) = \emptyset$. The associated 
Hopf co-Morita equivalent Hopf algebra is $GL_{-q}(2)$.
Indeed, one can show in the same way $GL_q(2,-2)$ is a 
right $GL_{-q}(2)$-comodule algebra
and is a right $GL_{-q}(2)$-Galois extension.
Therefore $GL_q(2,-2)$ is a $GL_q(2)-GL_{-q}(2)$ bigalois 
extension (\cite{[Sc]}, 3.4) and $GL_q(2)$ and $GL_{-q}(2)$ are 
Hopf co-Morita equivalent by \cite{[Sc]}, 5.7. 
A method to produce $GL_{-q}(2)$ from the data
 $(GL_q(2),GL_q(2,-2))$ is given in \cite{[Sc]}. 

\begin{rema}
{\rm The Hopf algebras $SL_q(2,{\mathbb C})$ et $SL_{-q}(2,{\mathbb
C})$ ($q$ not a root of unity or $q=1$) (\cite{[K]}) are not 
Hopf co-Morita equivalent.
Let us explain this result very briefly.

Let $\eta$ be a  fibre functor on Co$_f(SL_q(2,{\mathbb C}))$ such that
${\rm End}^{\vee}(\eta)
\widetilde {\longrightarrow} SL_{-q}(2,{\mathbb C})$ and let $V_q$ be
the obvious 2-dimensional comodule associated to
 $SL_q(2,{\mathbb C})$. The decomposition of tensor products of
irreducible  $SL_q(2,{\mathbb C})$-comodules (\cite{[K], [Ba1]}) 
shows that 
$\eta(V_q)\widetilde {\longrightarrow} V_{-q}$. Let us denote by $R_q$
the fundamental Yang-Baxter operator of 
$SL_q(2,{\mathbb C}$) (\cite{[K]}, VII.7).
An eigenvalue argument shows that $\eta(R_q) = -R_{-q}$ (with an obvious abuse of notation). It follows that we can find relations in
Hom$^\vee(\eta, \omega)$: there are generators 
$ (z_{ij})_{1\leq i,j\leq 2}$ satisfying the relations of $GL_q(2,-2)$, 
with the exception that the element $z_{11}z_{22} + q^{-1}z_{12}z_{21}$ 
must be a non-zero constant. But this element is not central, hence  
Hom$^\vee(\eta, \omega) = 0$, which contradicts the existence of the
fibre functor $\eta$.}
\end{rema}

\section{Unitary comodules}

In the remainder of the paper we assume the base field to be the field of complex numbers.

\subsection{Hopf $*$-algebras and conjugate comodules} 

\begin{defi}
A  {\rm Hopf $*$-algebra} is a Hopf algebra $A$ 
which is  a $*$-algebra in the usual sense and such that the comultiplication
$\Delta : A \longrightarrow A \otimes A$ is a $*$-homomorphism.
\end{defi}

It is well known that $\varepsilon$ is a $*$-homomorphism and that 
$S \circ * \circ S \circ *  = id$.

\medskip

\noindent
{\bf Notations}.
Let $V$ be a  complex vector space. The conjugate vector space of $V$
will be denoted by $\overline V$.
The canonical semi-linear isomorphism $V \longrightarrow
\overline V$ will be denoted by $j_V$. The canonical identification
 $V \longrightarrow \overline {\overline V}$ will be denoted by $\mu_V$, with
 $\mu_V = j_{\overline V} \circ j_V$.

\begin{defi}
Let $A$ be a Hopf $*$-algebra and let $V$ be an  $A$-comodule with
coaction $\alpha_V : V \longrightarrow V \otimes A$.
The {\rm conjugate comodule} of $V$ is the vector space
$\overline V$ endowed with the
 coaction $\alpha_{\overline V} : \overline V \longrightarrow
\overline V \otimes A$ defined by $\alpha_{\overline V} = (j_V \otimes *) \circ
\alpha_V \circ j_V^{-1} : \overline V \longrightarrow \overline V \otimes A$.
\end{defi}

The canonical identification $\mu_V : V \longrightarrow \overline {\overline V}$ is an isomorphism of comodules.
For every $f \in {\rm Hom}_{{\rm Co}_f(A)}(V,W)$, the conjugate map $\overline f = j_W \circ f
\circ j_V^{-1} : \overline V \longrightarrow \overline W$ is also
a map of comodules. Thus conjugation induces a semi-linear endofunctor
on ${\rm Co}_f(A)$.
Finally, let us note that for all comodules $V$ et $W$,
there is a natural monoidal isomorphism 
$\gamma_{V,W} : \overline{V \otimes W} \longrightarrow \overline W
\otimes \overline V$, $\overline {v \otimes w} \longmapsto \overline w \otimes
\overline v$.

\subsection{$*$-Categories}

A finite-dimensional Hilbert space is a pair
$(V,\phi_V)$ where
$V$ is a finite-dimensional vector space, and $\phi_V : \overline V \otimes
V \longrightarrow \mathbb C$ is a scalar product. Let us remark that in our conventions, the scalar product is linear in the second variable. By abuse if notation, we sometimes write $V$ for the pair $(V,\phi_V)$.

The category ${\cal H}$ of finite-dimensional Hilbert spaces is the
category whose objects are finite-dimensional
Hilbert spaces and whose maps are linear maps.

\begin{defi}
(\cite{[GLR]}). A {\rm $*$-category} is a 
$\mathbb C$-linear category  
$\mathcal C$, endowed with a semi-linear involutive contravariant 
functor $* : \mathcal C \rightarrow
\mathcal C$, which preserves  objects and satisfies 
the following two conditions:

1)
Let $X$ and $Y$ be  objects of $\mathcal C$ and let $f \in {\rm
Hom}_{\mathcal C}(X,Y)$. There is an element  $g \in {\rm Hom}_{\mathcal C}(X,X)$ such that $f^* \circ f = g^* \circ g$.

2) For every morphism $f$ of $\mathcal C$, the condition $f^* \circ f = 0$
implies $f = 0$.
\end{defi}

\begin{ex}
{\rm Endowed with the usual adjoint the category ${\cal H}$ is a 
$*$-category.}
\end{ex}

\begin{defi}
Let $\mathcal C$ and $\mathcal C'$  be $*$-categories. A {\rm $*$-functor}
 $F : \mathcal C
\longrightarrow \mathcal C'$  is a linear functor  such that 
$F(f^*) = F(f)^*$ for every morphism $f$ of $\mathcal C$.
\end{defi}

\subsection{Monoidal $*$-categories and unitary comodules}

Let $(V,\phi_V)$ and $(W,\phi_W)$ be finite-dimensional Hilbert spaces. 
Their Hilbert space tensor product is the Hilbert whose underlying vector space is $V \otimes W$ and whose scalar product $\phi_{V \otimes W}$ is given by
 $$\phi_{V \otimes W} = \phi_W \circ
(1_{\overline W}
\otimes \phi_V \otimes 1_W) \circ (\gamma_{V,W} \otimes 1_V \otimes 1_W)$$
(the map $\gamma_{V,W}$ is defined at the end of 3.1). Endowed with
this tensor product
${\cal H}$ is a monoidal category.

\begin{defi}
A {\sl monoidal $*$-category} is a $*$-category
$\mathcal C$, which is a monoidal category and such that for
any morphisms
$f$ and $g$ of $\mathcal C$, we have $(f \otimes g)^* = f^* \otimes g^*$.
\end{defi}

\begin{ex}
{\rm The category ${\cal H}$ is a monoidal $*$-category.}
\end{ex}
 
Let $(V,\phi_V)$ be a finite dimensional Hilbert space. Let $\kappa_V : \mathbb
C \longrightarrow V \otimes \overline V$ be the linear map defined by 
 $\kappa_V(1) = \sum_i e_i \otimes
\overline e_i$ where $(e_i)$ is an orthonormal basis of $V$. Then
$$(1_V \otimes \phi_V)\circ (\kappa_V \otimes 1_V) = 1_V \quad {\rm and}
\quad (\phi_V
\otimes 1_{\overline V}) \circ (1_{\overline V} \otimes \kappa_V) = 
1_{\overline V}.$$
This means that the triplet $(V,\phi_V,\kappa_V)$ is a left dual for $V$ in
${\cal H}$ (\cite{[JS]}, sec. 9). Therefore ${\cal H}$ is a 
rigid monoidal category.

\begin{defi}
Let $A$ be a Hopf $*$-algebra. A (finite-dimensional)
 {\rm unitary $A$-comodule} 
is a Hilbert space $(V,\phi_V)$ whose underlying vector space $V$ is an 
$A$-comodule and such that $\phi_V : \overline V \otimes V \longrightarrow
 \mathbb C$
is a map of comodules.

The {\rm category of unitary $A$-comodules}, denoted ${\cal U}_f(A)$,
is the category whose objects are unitary  $A$-comodules and   
whose maps are maps of comodules.
\end{defi}

An obvious monoidal structure is defined on ${\cal U}_f(A)$. Let 
$(V,\phi_V)$ be an unitary $A$-comodule.
The linear map $\kappa_V :
\mathbb C \longrightarrow V \otimes \overline V$ defined before
is a map of comodules. Hence if $f : V \rightarrow W$ is a map of unitary
comodules, then $f^* = (1_V \otimes \phi_W) \circ (1_V \otimes \overline
f \otimes 1_W) \circ (\kappa_V \otimes 1_W)$ is also a map of comodules.

Conclusion: ${\cal U}_f(A)$ is a rigid monoidal $*$-category.

\begin{defi}
Let $A$ be a Hopf $*$-algebra. A finite-dimensional
 $A$-comodule $V$ is said to be 
{\rm unitarizable} if there is a scalar product $\phi_V$ on $V$ such that
$(V,\phi_V)$ is a unitary $A$-comodule.
\end{defi}

\begin{defi}
A Hopf $*$-algebra $A$ is said to be {\rm unitarizable}
if every finite-dimensional $A$-comodule is unitarizable.
\end{defi}

\begin{ex}
{\rm The Hopf $*$-algebra of representative functions on
a compact quantum group is unitarizable (see 5.3.1).
Unitarizable Hopf $*$-algebras are called CQG algebras in \cite{[DK]}.
A Hopf $*$-algebra is unitarizable if and only if there is a positive
and faithful Haar measure on it (\cite{[DK]}, 3.10, in fact faithfulness of the Haar measure can be deduced from the other axioms: see \cite{[V2]}).}
\end{ex}

Let us remark that if $A$ is a unitarizable Hopf $*$-algebra, the inclusion 
functor ${\cal U}_f(A)
\subset {\rm Co}_f(A)$ is an equivalence of categories.

\subsection{$*$-Fibre functors}

\begin{defi}
Let $A$ be a unitarizable Hopf $*$-algebra. 
Let $\eta$ be a monoidal $*$-functor $ \mathcal
U_f(A)\longrightarrow \mathcal H$.
The functor $\eta$ is said to be a {\rm $*$-fibre functor
} on $\mathcal U_f(A)$ if the following isomorphisms (the monoidal constraints of $\eta$) are unitary isomorphisms:
$$\widetilde \eta_1 : {\mathbb C}\longrightarrow \eta({\mathbb C}
) \quad
 ; \quad \widetilde \eta_{V,W}: \eta(V)\otimes \eta(W)\longrightarrow \eta(V\otimes W).$$
\end{defi}

\begin{ex}
{\rm The forgetful functor $\omega : \mathcal U_f(A)\longrightarrow
\mathcal H$ is a  $*$-fibre functor.}
\end{ex}

Let us show that a  $*$-fibre functor is a particular fibre functor. 
The categories  $\mathcal U_f(A)$ and ${\rm Co}_f(A)$ are monoidally equivalent
 and thus a $*$-fibre functor induces
a monoidal functor (still denoted $\eta$) $\eta :{\rm Co}_f(A)\longrightarrow
{\rm Vect}_f({\mathbb C})$. The theory of orthogonal complements
shows that ${\rm Co}_f(A)$ is semisimple. Hence every $A$-comodule is a
direct sum of simple comodules (comodules $V$ such that End($V)= {\mathbb
C}$). The functor $\eta$ is linear and therefore it preserves  direct sums
and is faithful and exact.

\begin{defi}
Let $A$ be a unitarizable Hopf $*$-algebra. Let 
 $\eta$ and $\eta'$  be $*$-fibre functors on $\mathcal U_f(A)$.
A {\sl monoidal unitary isomorphism} $u : \eta \rightarrow \eta'$ is
a monoidal morphism between $\eta$ and $\eta'$ such that for every object
$V$ of $\mathcal U_f(A)$, the map $u_V:\eta(V)\rightarrow \eta'(V)$ is
unitary. The set of monoidal unitary isomorphisms between $\eta$ and $\eta'$
will be denoted by $U^\otimes(\eta,\eta')$.
The {\rm category of $*$-fibre functors on $\mathcal U_f(A)$}, denoted
 {\rm Fib$^*(A)$}, is the category whose objects are $*$-fibre functors on $\mathcal U_f(A)$ and whose morphisms are monoidal unitary isomorphisms.
\end{defi}

\subsection{An important Remark}

We end this section by an important remark from \cite{[W2]}, which will be useful in the next section.

  Let $(V,\phi_V)$ and $(V',\phi_{V'})$ be Hilbert spaces
and let $k : V \rightarrow V'$ be a semi-linear bijective map. 
Let $\alpha = k \circ
j_V^{-1} : \overline V \rightarrow V'$ be the associated linear isomorphism.
Let us define two linear maps (\cite{[W2]}, p. 39) 
$\overline t_k : V' \otimes V
\rightarrow \mathbb C$ and $t_k : \mathbb C \rightarrow V \otimes V'$ by
$\overline t_k = \phi_V \circ (\alpha^{-1} \otimes 1_V)$ and
 $t_k = (1_V \otimes
\alpha) \circ \kappa_V$. A direct computation shows that $t_{k^{-1}} =
(\overline t_k)^*$ and
$\overline t_{k^{-1}} = (t_k)^*$.

\section{Galois extensions and $*$-fibre functors}

 Let $A$ be a unitarizable Hopf $*$-algebra. Theorem 4.3.4 states the equivalence of categories between  $*$-fibre functors on $\mathcal U_f(A)$ and
$A$-$*$-Galois extension endowed with a positive Haar measure (see the definitions given in 4.1 and 4.2). The key results to construct the quasi-inverse functors are proposition 4.3.1 and 4.3.3. The end of the section is devoted to 
spectrum of an $A$-$*$-Galois extensions (4.4) and to Hopf $*$-co-Morita equivalence (4.5).

\subsection{$*$-Galois extensions}

\begin{defi}
Let $A$ be a Hopf $*$-algebra and let $Z$ be
an  $A$-Galois extension which is $*$-algebra and whose coaction  $\alpha_Z : Z \rightarrow A \otimes Z$ is a $*$-homomorphism. Then $Z$ is said to be 
an {\rm $A$-$*$-Galois extension}. 
\end{defi}

We will see in the proof of proposition 4.3.1
that if $Z$ is an $A$-$*$-Galois extension, the  associated fibre functor
(2.2) preserves conjugation.

\subsection{The Haar measure}

\begin{defi}
Let $A$ be a Hopf algebra over a field $k$
and let $Z$ be an
$A$-Galois extension. A {\rm Haar measure} on $Z$ is a linear map  $\mu :
Z \rightarrow k$ such that
$u_A \circ \mu = (1_A \otimes \mu)\circ \alpha_Z$ and $\mu(1_Z) = 1_k$.
\end{defi}

If $Z = A$, a Haar measure on $A$ is a Haar measure in the usual sense
(\cite{[A]}). A Haar measure $J$ on $A$ is unique and right invariant: 
$(J \otimes 1_A) \circ \Delta =u_A \circ J$. There is a Haar measure 
on $A$ if and only if the category Co$_f(A)$ is semisimple.

\begin{prop}
Let $A$ be a Hopf algebra endowed with a Haar measure and let $Z$ be
an  $A$-Galois extension. There is a unique Haar measure on $Z$.
\end{prop}

\noindent
{\bf Proof}. Existence. Let $f$ be a linear form on $A$ such that $f(1_Z)=1_k$ and let $J$ be a Haar measure on $A$. An easy computation shows that $J*f= (J\otimes f) \circ \alpha_Z$ 
is a Haar measure on $Z$.

\noindent
Uniqueness.  Let $J$ be a Haar measure on $A$. The right invariance of $J$ implies that 
 $\alpha_Z \circ ((J \otimes 1_Z) \circ \alpha_Z) =u_A \otimes 1_Z \circ ((J \otimes 1_Z)\circ \alpha_Z)$. Therefore for all $z \in Z$, we have
$(J \otimes 1_Z)\circ \alpha_Z(z) \in k$ since $k \wedge Z \cong k$.
Thus if $\mu$ and $\mu'$ are Haar measures on $Z$, we have
 $\mu = J*\mu =J*\mu'=\mu'$. \square   

\begin{defi}
Let $A$ be a Hopf $*$-algebra and let $Z$ be an $A$-$*$-Galois extension.
A Haar measure $\mu$ on $Z$ said to be {\rm positive} if
$\mu(z^*z) \geq  0$ for all $z$. It is said to be {\rm faithful} if 
$\mu(z^*z) >0 $ for all  $z \not = 0$.
\end{defi}

Let $A$ be a Hopf $*$-algebra. There is a positive Haar measure on $A$
if and only if $A$ is unitarizable (combine \cite{[DK]}, 3.10 and \cite{[V2]}).

\begin{defi}
Let $A$ be a unitarizable Hopf $*$-algebra. 
The {\rm category Gal$^*(A)$} is the category whose objects are
$A$-$*$-Galois extension with a positive Haar measure and
whose morphisms are $*$-homomorphisms which are morphisms of $A$-Galois 
extensions.
\end{defi}

The following result shows that the positivity condition on the Haar measure 
is not restrictive when we have a view towards compact quantum groups. 
 
\begin{prop}
 Let $A$ be a unitarizable Hopf $*$-algebra and let
$Z$  be an $A$-$*$-Galois extension. Assume that there is a positive linear 
form on $Z$  ($\psi(z^*z)\geq 0$ for all $z$) such that
$\psi(1)=1$. Then there is a positive Haar measure on $Z$.
\end{prop}

\noindent 
{\bf Proof}. Let $J$ be a  positive and faithful Haar measure on $A$ 
(\cite{[DK]}, 3.10) and let
$\mu = J* \psi = (J\otimes \psi) \circ \alpha_Z$ be the Haar measure on
$Z$. Let $z \in Z$ with $\alpha_Z(z) = \sum a_i\otimes t_i$. Endowed with
the scalar product
$\langle a, b\rangle = J(a^*b)$, then  $A$ is a prehilbert space.
One can assume that the  $a_i$'s are  orthonormals. Hence $\mu(z^*z) = \sum
\psi(t_i^*t_i) \geq 0$  and  $\mu$ is a positive Haar measure
on $Z$. \square

\subsection{The equivalence of categories}

In the next result, Ulbrich's functor (2.2) is adapted to our setting.

\begin{prop}
 Let $A$ be a unitarizable Hopf $*$-algebra and let
$Z$  be an $A$-$*$-Galois extension endowed with a positive Haar measure.
Let $\eta_Z : {\rm Co}_f(A) \longrightarrow {\rm Vect}_f({\mathbb
C)}$ be the associated fibre functor (2.2). Then $\eta_Z$ factorizes through
a $*$-fibre functor $\eta_Z^\odot : {\cal U}_f(A) \longrightarrow {\cal H}$
followed by the forgetful functor.
\end{prop}

\noindent
{\bf Proof}. {\bf Step 1}. Let $V$ be an $A$-comodule. We are going to
define an  isomorphism
$\lambda_V : \overline{V \wedge Z} \longrightarrow \overline V \wedge Z$. 
Let $\lambda_V \left ( \overline{\sum_i v_i \otimes z_i} \right ) = \sum_i
\overline v_i \otimes z^*_i$: 
we have $\lambda_V(\overline{V \wedge Z}) \subset
\overline V \wedge Z$ since $\alpha_Z : Z \longrightarrow A \otimes Z$ is a
$*$-homomorphism. It is easily seen that $\lambda_V$ is an isomorphism.

\smallskip
\noindent 
{\bf Step 2}. Let $(V,\phi_V)$ be a unitary $A$-comodule. We are going 
to define a 
scalar product on $V \wedge Z$. Let $\phi'_V$ be the linear map:

$\overline{V \wedge Z} \otimes V \wedge Z \stackrel{\lambda_V \otimes
1}{\longrightarrow} \overline V \wedge Z \otimes V \wedge Z {\widetilde
{\longrightarrow}} (\overline V \otimes V) \wedge Z \stackrel{\phi_V \otimes
1_Z}{\longrightarrow} {\mathbb C} \wedge Z \tilde \rightarrow {\mathbb C}$.

\noindent
Let us show that $\phi'_V$ is a scalar product on $V \wedge Z$.
Choose $(v_i)_{1\leq i \leq n}$ an orthonormal basis of $V$ (with respect to
$\phi_V$) and let $z = \sum_i v_i \otimes z_i$ and $z' = \sum_i v_i \otimes z'_i$ be two elements of $V \wedge Z$. Then $\phi'_V(\overline z \otimes z') =
\alpha$ where $(\sum_i z^*_i z'_i) =  \alpha 1_Z$ with $\alpha \in \mathbb C
$ since $\mathbb C\wedge Z {\widetilde {\longrightarrow}} \mathbb C$. 
Thus if $\mu$ denotes the Haar measure on $Z$, we have $\phi'_V(\overline z
\otimes z') = \mu(\sum_i z^*_i z'_i)$, and hence  
$\phi'_V(\overline z
\otimes z') = \overline {\phi'_V (\overline z' \otimes z)}$ and
$\phi'_V(\overline z \otimes z) \geq  0$ since $\mu$ is positive.
The linear map $\overline {V\wedge Z} \rightarrow(V\wedge Z)^*$ induced by
$\phi'_V$ is bijective since it is the composition 
of the following isomorphisms:  

$\overline{V \wedge Z} \stackrel{\lambda_V
}{\longrightarrow} \overline V \wedge Z  {\widetilde
{\longrightarrow}} ( V^* \wedge Z) {\widetilde{\longrightarrow}} (V \wedge Z)^*$ (the last map is given by the unicity up to isomorphism of the dual in a
monoidal category).

\noindent 
It follows that $\phi'_V$ is a scalar product since it is sesquilinear, 
positive and non-degenerate.

A functor $\eta^\odot_Z : {\cal U}_f(A) \longrightarrow {\cal H}$
is thus defined in the following way :
$\eta^\odot_Z(V) = (V \wedge Z, \phi'_V)$ and $\eta^\odot_Z(f) = f \otimes 1_Z$.

\smallskip
\noindent
{\bf Step 3}. Let us show now that $\eta^\odot_Z$ is a $*$-functor. 
Let $(V,\phi_V)$ and $(W,\phi_W)$ be unitary $A$-comodules and let $f : V
\longrightarrow W$ be a map of
$A$-comodules. Let us recall that $\eta^\odot_Z(f) = f \otimes 1_Z$ and that
$\eta^\odot_Z(f)^*$ is the only linear map such that 
$$\phi'_W \circ
(\overline{\eta_Z^\odot(f)} \otimes 1_{\eta_Z(W)}) = \phi'_V \circ (1_{\eta_Z(V)}
\otimes \eta_Z^\odot(f)^*).$$
We have $\overline{\eta_Z^\odot(f)} =
\lambda^{-1}_W \circ \eta_Z^\odot(\overline f) \circ \lambda_V$ ($\lambda_V$
is defined in step 1), and an easy computation now shows that
$\eta^\odot_Z(f^*) = \eta^\odot_Z(f)^*$.

\smallskip
\noindent
{\bf Step 4}. It remains to prove that the isomorphisms $\mathbb C {\widetilde
{\longrightarrow}} \mathbb C \wedge Z$ et $(V \wedge Z) \otimes (W \wedge Z)
{\widetilde {\longrightarrow}} (V \otimes W) \wedge Z$ are unitaries. For
the first one, it is obvious. Let $(V,\phi_V)$ and $(W,\phi_W)$ be unitary 
$A$-comodules
with orthonormal bases $(v_i)_{1 \leq i \leq n}$ and $(w_i)_{1 \leq i
\leq n'}$ respectively. Let $(\sigma_s)_{1 \leq s \leq p}$ and $(\tau_r)_{1
\leq r \leq p'}$ be orthonormal bases of $V \wedge Z$ and $W \wedge Z$
respectively: $\sigma_s =\sum_i v_i \otimes z_{is}$ and $\tau_r = \sum_i
w_i \otimes t_{ir}$. We have $\sum_i z^*_{is} z_{is'} = \delta_{ss'}$ and
$\sum_i t^*_{ir} t_{ir'} = \delta_{rr'}$. It is then easily seen that
the isomorphism $(V \wedge Z) \otimes (W \wedge Z) \longrightarrow (V \otimes W)
\wedge Z$ preserves  orthonormal bases and therefore is  unitary.\square 

\bigskip

The following lemma is useful to describe the Haar measure on a Galois extension.

\begin{lemm}
Let $A$ be a Hopf algebra over an algebraically closed field $k$ and 
let $\eta : {\rm Co}_f(A) \longrightarrow {\rm Vect}_f(k)$ be fibre functor. 
Assume that there is a Haar measure on $A$ (${\rm Co}_f(A)$ is semisimple). 
Let us denote by $\widehat A$ the set of isomorphism classes of irreducible 
$A$-comodules and choose for every $\lambda \in \widehat A$ a 
representative $V_{\lambda}$. Then there is a linear 
isomorphism $I : \bigoplus_{\lambda \in \hat A} {\rm
Hom}(\omega(V_\lambda), \eta(V_\lambda)) {\widetilde {\longrightarrow}} {\rm
Hom}^\vee(\eta,\omega)$.
\end{lemm}

\noindent
{\bf Proof}. The map $I$ is defined as follows: 
if $v \in {\rm Hom}(\omega(V_\lambda),\eta(V_\lambda))$, let
$I(v) = [V_\lambda, v]$. If $V$ is an object of ${\rm Co}_f(A)$, there are
elements $\lambda_1, \ldots, \lambda_n$ of $\widehat A$ such that $V \cong
\bigoplus^n_{i=1} V_{\lambda_i}$. Let $u_i : V_{\lambda_i} \longrightarrow V$
and $p_i : V \longrightarrow V_{\lambda_i}$ be the arrows of direct sum. Let
$v \in {\rm Hom}(\omega(V),\eta(V))$, then $[V,v] = 
[V,\eta(\sum_i u_i \circ p_i) \circ
v] = \sum_i [V_{\lambda_i}, \eta(p_i) \circ v \circ u_i]$, therefore $I$ is
surjective.

Let $\lambda \in \widehat A$ and let $v \in {\rm
Hom}(\eta(V_\lambda),\omega(V_\lambda))$. 
It is easy to see that there is an element 
$u$ from ${\rm Hom}(\eta,\omega)$ such that $u_{V_\lambda} = v$ and  $u_{V_\mu} = 0$ if $\mu \not = \lambda$. Let $x = \sum_j [V_{\lambda_j}, u_j] \in
{\rm Hom}^\vee(\eta,\omega)$. the preceding assertion
shows that for every 
$j$ and $v \in {\rm Hom}(\eta(V_{\lambda_j}),\omega(V_{\lambda_j}))$, there is a
linear form $f$ (see 2.5) such that $f(x) =  {\rm Tr}(v \circ u_j)$. If $x =0$,
then $u_j = 0$ for all $j$ and therefore $I$ is injective. \square

\begin{prop}
Let $A$ be a unitarizable  Hopf $*$-algebra and let
$\eta : {\cal U}_f(A) \longrightarrow {\cal H}$ be a  $*$-fibre functor. Let
$Z$ be the $A$-Galois extension associated to the fibre functor on 
${\rm Co}_f(A)$ defined by $\eta$ 
(2.3). Then $Z$ is an $A$-$*$-Galois-extension 
endowed with a positive and faithful Haar measure.
\end{prop}

\noindent
{\bf Proof}. By 2.3 we have $Z = {\rm Hom}^\vee(\eta,\omega)$ 
and the Hopf algebras  $A$ and ${\rm End}^\vee (\omega)$ (2.3) are identified 
($\omega$ is the forgetful functor).

\smallskip
\noindent
{\bf Step 1}. Let us describe the Haar measure on ${\rm
Hom}^\vee(\eta,\omega)$. We know that ${\rm Co}_f(A)$ is semisimple.
Let $\widehat A$ be the set of isomorphism classes of irreducible comodules.
For every $\lambda \in \widehat A$, let us choose a representative 
$V_\lambda$. 
By lemma 4.3.2 we can define  $\mu = \bigoplus_{\lambda \in \hat A} \mu_\lambda$ where $\mu_\lambda =
0$ if $\lambda \not = 1$ (the trivial comodule class) and $\mu_1 = 1$.
It is easily seen that $\mu$ is the Haar measure.

\smallskip
\noindent
{\bf Step 2}. We now construct the involution on ${\rm Hom}^\vee(\eta,\omega)$ . It is given by $[V,u]^* = [\overline V, \overline u]$
 on the Hopf algebra ${\rm End}^\vee(\omega)$.

Let $(V,\phi_V)$ be a unitary comodule. Let 
$$
E_V = \widetilde \eta^{-1}_1 \circ \eta(\phi_V) \circ \widetilde 
\eta_{\overline V, V} :
\eta(\overline V) \otimes \eta(V) \longrightarrow \mathbb C
$$
$$
K_V = \widetilde \eta^{-1}_{V,\overline V} \circ \eta(\kappa_V) \circ \widetilde
\eta_1 :
\mathbb C
\longrightarrow \eta(V) \otimes \eta(\overline V)
$$
(where $\kappa_V$ is the duality map defined in 3.3).
Let $\eta(V,\phi_V) = (\eta(V), \phi'_V)$. By the unicity of duals in
monoidal categories, there is a unique isomorphism $\lambda_V :
\overline{\eta(V)} \longrightarrow \eta(\overline V)$ such that $\phi'_V = E_V
\circ (\lambda_V \otimes 1_{\eta(V)})$ and $\kappa'_V = (1_{\eta(V)}
\otimes \lambda^{-1}_{V})
\circ K_V$ (with $\kappa'_V$ as in 3.3).

Let $f : (V,\phi_V) \longrightarrow (W,\phi_W)$ be a map of comodules,
Then we have $\lambda_W \circ \overline{\eta(f)} = \eta(\overline f) \circ \lambda_V$.
Indeed $\overline {\eta(f)}$ is the only map such that $\phi'_W \circ
(\overline{\eta(f)} \otimes 1_{\eta(W)}) = \phi'_W \circ (1_{\overline {\eta(V)}}
\otimes \eta(f)^*)$: it is easy to check that 
$\lambda_W^{-1} \circ \eta (\overline f) \circ \lambda_V$ satisfies this equation
(use $\eta(f^*) = \eta(f)^*$). 

This fact allows us to define a semi-linear endomorphism of
${\rm Hom}^\vee(\eta,\omega)$ 
denoted $\sigma : \sigma([V,v]) = [\overline V, \lambda_V \circ \overline v]$.
The isomorphisms $\lambda_V$ are monoidal in
an obvious way since they are unique.
Hence $\sigma$ is anti-multiplicative (we leave the formal verification to the reader).

It remains to show that $\sigma$ is involutive. 
Let $(V,\phi_V)$ be an unitary $A$-comodule. We have
 $\sigma^2([V,v]) =
[\overline {\overline V}, \lambda_{\overline V} \circ \overline \lambda_V \circ
\overline {\overline v}] = [V,
\eta(\mu_V^{-1}) \circ \lambda_{\overline V} \circ \overline \lambda_V \circ
\overline {\overline v} \circ \mu_V] = [V, \eta(\mu_V^{-1}) \circ
\lambda_{\overline V}
\circ \overline \lambda_V \circ \mu_{\eta(V)} \circ v]$.
Therefore it is sufficient to show that $\eta(\mu^{-1}_V) \circ
\lambda_{\overline V} \circ \overline \lambda_V \circ \mu_{\eta(V)} =
1_{\eta(V)}$. Let $k = \lambda_V \circ j_{\eta(V)} : \eta(V) \longrightarrow
\eta(\overline V)$: $k$ is a semi-linear isomorphism. We use 3.5 now.
We have:
$$
\overline t_k = \phi'_V \circ (\lambda_V^{-1} \otimes 1_{\eta(V)}) = E_V :
\eta(\overline V)
\otimes \eta(V) \longrightarrow \mathbb C
$$
$$
t_k = (1_{\eta(V)} \otimes \lambda_V) \circ \kappa'_V = K_V : \mathbb C
\longrightarrow
\eta(V) \otimes \eta(\overline V)
$$
(with the same notations as in the beginning of step 2).
Moreover
$$
\overline t_{k^{-1}} = \phi'_{\overline V} \circ ((j_{\eta(\overline V)} \circ
\lambda_V
\circ j_{\eta(V)}) \otimes 1_{\eta(\overline V)}) : \eta(V) \otimes
\eta(\overline V) \longrightarrow \mathbb C \leqno (+)
$$
$$
t_{k^{-1}} = (1_{\eta(\overline V)} \otimes (j^{-1}_{\eta(V)} \circ
\lambda_V^{-1} \circ j_{\eta(\overline V)}^{-1})) \circ \kappa'_{\overline V} :
\mathbb C
\longrightarrow \eta (\overline V) \otimes \eta(V).
$$
By 3.5 we have  $\overline t_{k^{-1}} = (t_k)^* = \widetilde
\eta_1^{-1} \circ
\eta (\kappa^*_V) \circ \widetilde \eta_{V,\overline V}$ (we use the fact that
$\widetilde
\eta_1$ et $\widetilde \eta_{V,\overline V}$ are unitaries).
It is immediate that $\kappa^*_V = \phi_{\overline V} \circ \mu_V \otimes 1_{\overline V}$
, hence

\begin{eqnarray*}
&  \overline t_{k^{-1}} &  = \widetilde \eta_1^{-1} \circ \eta(\phi_{\overline
V})
\circ \widetilde \eta_{\overline {\overline V}, \overline V} \circ (\eta(\mu_V)
\otimes 1_{\eta(\overline V)}) \\
 & & = E_{\overline V} \circ (\eta(\mu_V) \otimes
1_{\eta(\overline V)}) = \phi'_{\overline V} \circ ((\lambda_V^{-1} \circ
\eta(\mu_V)) \otimes 1_{\eta(\overline V)}).
\end{eqnarray*}
Combining this equality with (+), and using the fact that $\phi'_{\overline V}$is non degenerate, we get

$ j_{\eta(\overline V)} \circ \lambda_V \circ j_{\eta(V)}  =
\lambda^{-1}_{\overline V} \circ \eta(\mu_V)$,  
hence $\overline \lambda_V \circ j_{\overline {\eta(V)}} \circ
j_{\eta(V)} = \lambda^{-1}_{\overline V} \circ \eta(\mu_V)$. 

\noindent
This shows that
$
\eta(\mu^{-1}_V) \circ \lambda_{\overline V} \circ \overline \lambda_V
\circ \mu_{\eta(V)} = 1_{\eta(V)}$, and thus $\sigma$ turns
${\rm Hom}^\vee(\eta,\omega)$
into a $*$-algebra.

It is easy to see that the coaction ${\rm Hom}^\vee(\eta,\omega)
\longrightarrow {\rm End}^\vee(\omega) \otimes {\rm
Hom}^\vee(\eta,\omega)$ is $*$-homomorphism.

\medskip
\noindent
{\bf Step 3}. It remains to show that the Haar measure is positive and faithful. Let $a \in
{\rm Hom}^\vee(\eta,\omega)$ with $a = \sum_i a_i$ where $a_i = [V_i,v_i]$ 
and the $V_i$'s are distinct irreducible comodules (see lemma 4.3.2). 
If $i \not = j$ we have
$\mu(a^*_i a_j) = 0$ since $\mathbb C$ does not appear as a subcomodule of 
$\overline V_i \otimes V_j$ (By duality 
${\rm Hom}(\mathbb C,\overline V_i \otimes V_j) {\widetilde {\longrightarrow}} {\rm
Hom}(V_i,V_j)$).
Therefore it is sufficient to show  that $\mu(a^* a) > 0$
with $a = [V,v] \not = 0$ and $V$ irreducible.

We have $a^*a = [\overline V \otimes V, \widetilde \eta_{\overline V, V} \circ
((\lambda_V \circ \overline v) \otimes v)]$ ($\lambda_V$ is defined in step 2).

Let $\phi_V$ be an invariant scalar product on $V$
and let $\phi'_V$ be a scalar product on $\eta(V)$ such that
$\eta(V,\phi_V) = (\eta(V),\phi'_V)$. Let 
$\alpha = \phi_V \circ \phi^*_V \in \mathbb C$
:  $\alpha >0$ (if 
$A$ is commutative then $\alpha = \dim V$). Let $p = \alpha^{-1} \phi_V^*
\circ \phi_V$: $p$ is a projection in ${\rm End}(\overline V \otimes V)$.
Hence we have 
$$
a^*a = [\overline V \otimes V, \eta(p) \circ \widetilde \eta_{\overline V, V}
\circ ((\lambda_V \circ \overline v)\otimes v)] + [\overline V \otimes V,
\eta(1-p)
\circ
\ldots]
$$

The trivial comodule $\mathbb C$ appears exactly once in the decomposition
of $\overline V\otimes V$ into irreducible comodules. Hence
$\mu([\overline V\otimes V, \eta(1-p)\circ
\ldots])=0$. On the other hand

$ [\overline V \otimes V, \eta(p)\circ \widetilde \eta_{\overline V,V} \circ
((\lambda_V
\circ \overline v)
\otimes v)]
  = \alpha^{-1}[\mathbb C, \eta(\phi_V) \circ \widetilde \eta_{\overline
V, V} \circ ((\lambda_V \circ \overline v) \otimes v) \circ \phi^*_V]$ 

$  = \alpha^{-1} (\tilde \eta_1^{-1} \circ \eta(\phi_V) \circ \widetilde
\eta_{\overline V,V} \circ (\lambda_V \otimes 1_{\eta(V)}) \circ (\overline v
\otimes 1_{\eta(V)})
\circ (1_{\overline V} \otimes v) \circ \phi^*_V). [\mathbb C,\widetilde \eta_1]$

$= \alpha^{-1} (\phi'_V \circ (\overline v \otimes 1_{\eta(V)}) \circ
(1_{\overline V} \otimes v) \circ \phi^*_V) [\mathbb C, \widetilde \eta_1]$ 

$= \alpha^{-1} (\phi_V \circ (1_{\overline V} \otimes v^*) \circ (1_{\overline
V} \otimes v) \circ \phi^*_V) [\mathbb C, \widetilde \eta_1]$ 

$ =\alpha^{-1} ((1_{\overline V} \otimes v) \circ \phi^*_V)^* \circ
((1_{\overline V} \otimes v) \circ \phi^*_V) [\mathbb C, \widetilde \eta_1]$.

\noindent
This constant is  strictly positive (see definition 3.2.1), therefore the Haar measure is
 positive and faithful. \square

\bigskip

We are now ready to state our equivalence of categories.
Let $A$ be a  unitarizable Hopf $*$-algebra. Recall that ${\rm Gal}^*(A)$
is the category of $A$-$*$-Galois extensions endowed with a positive Haar 
measure (definition 4.2.4) and that ${\rm Fib}^*(A)$ the category of $*$-fibre 
functors on ${\cal U}_f(A)$ (definition 3.4.3).

\begin{theo}
Let $A$ be a unitarizable Hopf $*$-algebra. There is a functor
${\rm Gal}^*(A) \longrightarrow {\rm Fib}^*(A)$ (given on the objects by
proposition 4.3.1) which is an equivalence of categories.
\end{theo}

\noindent
{\bf Proof}. {\bf Step 1}. Let us first define the functor. Let $Z$ be 
an $A$-$*$-Galois extension endowed with a positive Haar measure
and let $\eta^\odot_Z : {\cal U}_f(A) \longrightarrow
{\cal H}$ the $*$-fibre functor constructed in proposition 4.3.1.
We use the functor from 2.2 : ${\rm Fib}(A) \longrightarrow {\rm Gal}(A)$.
Let  $f : Z \longrightarrow T$ be a morphism of ${\rm Gal}^*(A)$. 
A morphism of fibre functors $u : \eta^\odot_Z \longrightarrow \eta^\odot_T$ 
is associated to $f$ in the following way:
let $(V,\phi_V)$ be a unitary 
$A$-comodule, $u_V = 1_V \otimes f : V \wedge Z \longrightarrow V \wedge T$.
We must show that $u_V$ is unitary
(for the scalar products defined in step 2 of the proof of proposition 4.3.1).
Choose $(v_i)_{1 \leq i \leq n}$ an orthonormal basis of $V$ and $(\sigma_j)_{1
\leq j \leq p}$ an orthonormal basis of $V \wedge Z$. 
Let  $\sigma_j =
\sum^n_{i=1} v_i \otimes z_{ij}$ where $(z_{ij})_{1 \leq i \leq n, 1 \leq j \leq
p}$ are elements of $Z$. The basis $(\sigma_j)_{1 \leq j \leq p}$
is orthonormal if and only if
$\sum_i z^*_{ij} z_{ik} = \delta_{jk}$ for all $j$ and
$k$. Since $f$ is a $*$-homomorphism the map $u_V$  preserves orthonormal bases and 
is unitary.

Thus we have defined a functor ${\rm Gal}^*(A) \longrightarrow {\rm
Fib}^*(A)$.

\medskip
\noindent
{\bf Step 2}. Let us now define a functor  ${\rm Fib}^*(A) \longrightarrow {\rm Gal}^*(A)$. We use the functor defined in 2.3 and proposition 4.3.3. 
Let $\eta$ and $\theta$ be $*$-fibre functors, and let $u
\in U^\otimes(\eta,\theta)$ a monoidal unitary isomorphism.
We must show that the morphism of Galois extension
$$
f : {\rm Hom}^\vee(\eta,\omega) \longrightarrow {\rm Hom}^\vee(\theta,\omega), \
[V,v]
\longmapsto [V, u_V \circ v]
$$
is a $*$-homomorphism for the involutions constructed of step 2 in the proof of proposition 4.3.3.
We use the notations of step 2 in the proof of proposition 4.3.3 with a mark
$\eta$ or $\theta$ to label the functor.

Let $(V,\phi_V)$ be a unitary $A$-comodule with $\eta(V,\phi_V) =
(\eta(V),\phi^{\eta'}_V)$ and $\theta(V,\phi_V) = (\theta(V),
\phi^{\theta'}_V)$.
We have $f([V,v]^*) = f[\overline V,\lambda^\eta_V \circ \overline v]) = [\overline
V,u_{\overline V} \circ \lambda^\eta_V \circ \overline v]$ and $f([V,v])^* =
[V,u_V \circ v]^* = [\overline V, \lambda^\theta_V \circ \overline u_V \circ
\overline v]$. Therefore it is enough to see
that $u_{\overline V} \circ \lambda_V^\eta = \lambda_V^\theta \circ \overline
u_V$. The following equalities hold :
$$u_{\mathbb C} \circ \widetilde \eta_1 = \widetilde \theta_1 \quad ;\quad u_{\mathbb C} \circ
\eta(\phi_V) =
\theta(\phi_V) \circ u_{\overline V \otimes V} 
\quad ;\quad u_{\overline V \otimes V}
\circ \widetilde \eta_{\overline V,V} = \widetilde \theta_{\overline V,V} \circ
(u_{\overline V} \otimes u_V)$$ 
 $$\phi^{\theta'}_V = \widetilde \theta^{-1}_1 \circ
\theta(\phi_V) \circ \widetilde \theta_{\overline V,V} \circ (\lambda_V^\theta
\otimes 1_{\theta(V)})\quad ; \quad 
\phi^{\eta'}_V = \widetilde \eta_1^{-1} \circ
\eta(\phi_V)
\circ \widetilde \eta_{\overline V,V} \circ (\lambda^\eta_V \otimes
1_{\eta(V)})$$
Since $u$ is unitary, we have $\phi_V^{\theta'} \circ \overline u_V \otimes u_V =
\phi^{\eta'}_V$. This equality implies:
\begin{eqnarray*}
\theta(\phi_V)  \circ \widetilde \theta_{\overline V, V} \circ
(\lambda^\theta_V \otimes 1_{\eta(V)}) \circ (\overline u_V \otimes u_V) &  = &
u_{\mathbb C} \circ \eta(\phi_V) \circ \widetilde \eta_{\overline V,V} \circ
(\lambda^\eta_V
\otimes 1_{\eta(V)}) \\
&  = & \theta(\phi_V) \circ u_{\overline V \otimes V} \circ \widetilde
\eta_{\overline V,V} \circ (\lambda^\eta_V \otimes 1_{\eta(V)}) \\
& = & \theta(\phi_V) \circ \widetilde \theta_{\overline V,V} \circ (u_{\overline
V}
\otimes u_V)
\circ (\lambda^\eta_V \otimes 1_{\eta(V)}) \\
 {\rm hence} \quad  \theta(\phi_V) \circ \widetilde \theta_{\overline V,V} \circ
((\lambda^\theta_V \circ \overline u_V) \otimes 1_{\eta(V)}) & =  & \theta(\phi_V)
\circ \widetilde \theta_{\overline V,V} \circ ((u_{\overline V} \circ
\lambda^\eta_V)
\otimes 1_{\eta(V)}).
\end{eqnarray*}
The map $\theta(\phi_V) \circ \tilde \theta_{\overline V,V}$ is non-degenerate
and therefore we have the desired equality.

Thus the functor from 2.3 ${\rm Fib}(A) \longrightarrow {\rm Gal}(A)$ 
gives rise to a functor  ${\rm Fib}^*(A) \longrightarrow {\rm Gal}^*(A)$.

\medskip
\noindent
{\bf Step 3}. Let us show that the functors defined in step 1 and step 2 are
 quasi-inverse. We use 2.4. 

Let $Z$ be $A$-$*$-Galois extension endowed with a positive Haar measure. Let
$\eta^\odot_Z$ be the associated $*$-fibre functor and let $f : {\rm
Hom}^\vee(\eta_Z^\odot,\omega) \longrightarrow Z$ be the isomorphism of 2.4.
Let $(V,\phi_V)$ be an unitary $A$-comodule and let $x \in V \wedge Z$ and $\phi \in
V^*$, we have $f([V,\phi \otimes x]) = (\phi \otimes 1_Z) \circ (1_V \otimes
\varepsilon \otimes 1_Z) \circ (\alpha_V \otimes 1_Z)(x)$. 
The reader will easily check,
using orthonormal bases, that $f$ is a
$*$-homomorphism (the isomorphism $\overline {V \wedge Z} \longrightarrow
\overline V \wedge Z$ is given in step 1 of the proof of proposition 4.3.1).

\medskip
Let $\theta$ be a $*$-fibre functor. Let us recall how the fibre functor
isomorphism $\gamma : \theta \longrightarrow \eta^\odot_{{\rm
Hom}^\vee(\theta,\omega)}$ from 2.4 is constructed.
Let $(V,\phi)$ be a unitary comodule with orthonormal basis 
$(v_i)_{1\leq i \leq n}$ and let $(x_j)_{1 \leq i \leq p}$ be an orthonormal
basis of $(\theta(V),\phi'_V)$. We have
$\gamma_V(x_j) = \sum_i v_i \otimes [V, v^*_i \otimes x_j]$. Let
$z_{ij} = [V, v^*_i \otimes x_j]$. To show that $\theta$ is unitary, it is
enough to see that $\sum_i z_{ij}^* z_{ik} = \delta_{jk}$ for all $j$ and $k$
(see step 2 in the proof of proposition 4.3.1).
We have $z^*_{ij} z_{ik} = [\overline V \otimes V, \widetilde \theta_{\overline V,V}
\circ ((\lambda_V \circ (\overline{v_i^* \otimes x_j})) \otimes (v^*_i \otimes
x_k)] = [\overline V \otimes V, \widetilde \theta_{\overline V,V} \circ (\lambda_V
\otimes 1_{\eta(V)}) \circ ((\overline v_i^* \otimes v_i^*) \otimes (\overline x_j
\otimes x_k))]$. Hence
$$
\sum^n_{i=1} z^*_{ij} z_{ik} = [\overline V \otimes V, \widetilde
\theta_{\overline V,V} \circ (\lambda_V \otimes 1_{\eta(V)}) \circ (\overline x_j
\otimes x_k) \circ
\phi_V] = [\mathbb C, \widetilde \theta_1\circ \phi'_V \circ \overline x_j \otimes
x_k] = $$
$\delta_{jk} [\mathbb C, \tilde \theta_1] 
= \delta_{jk} 1$.
The proof is now complete.  \square

\begin{coro}
Let $A$ be a unitarizable Hopf $*$-algebra and let $Z$ be an $A$-$*$-Galois
extension. A positive Haar measure on $Z$ is faithful.
\end{coro}

\noindent
{\bf Proof}. Assume there is a positive Haar measure on $Z$. By theorem 4.3.4 and proposition 4.3.3 $Z$ is isomorphic with an $A$-$*$-Galois
extension whose positive Haar measure is faithful. The Haar measure on 
$Z$ is unique and hence must be faithful. \square

\subsection{The spectrum of a $A$-$*$-Galois extension}

\begin{theo}
Let $A$ be a unitarizable Hopf $*$-algebra and let $\eta
: {\cal U}_f(A) \longrightarrow {\cal H}$ be a  $*$-fibre functor. Then there
is a bijection ($\omega$  is the forgetful functor):
$$
U^\otimes (\eta,\omega) \longrightarrow {\rm Hom}_{*-{\rm alg}} ({\rm
Hom}^\vee(\eta,\omega), \mathbb C).
$$
\end{theo}

\noindent
{\bf Proof}. Recall that $\eta$ may be seen as fibre functor on 
${\rm Co}_f(A)$. Hence by 2.5 there is a bijection :
$$
{\rm Hom}^\otimes (\eta,\omega) \longrightarrow {\rm Hom}_{\mathbb C-{\rm alg}} ({\rm
Hom}^\vee(\eta,\omega), \mathbb C)
$$
which associates to $u \in {\rm Hom}^\otimes(\eta,\omega)$ the linear form $f_u$
with $f_u (\sum  [V_i,v_i]) = \sum_i {\rm Tr}(u_{V_i} \circ v_i)$.
We will use the notations of step 2 in the proof of proposition 4.3.1 (and of 
steps  2 and 3 in the proof of theorem 4.3.4).

Let $u \in U^\otimes (\eta,\omega)$ and let $(V,\phi_V)$ be a unitary $A$-comodule. We have
 $f_u ([V,v]^*) = {\rm Tr}(u_{\overline V} \circ \lambda_V \circ
\overline v)$ and
$\overline {f_u([V,v])} = \overline {{\rm Tr}(u_V \circ v)}$. We want to show
that $f_u$
is a $*$-homomorphism. It is enough to see that
$u_{\overline V} \circ \lambda_V = \overline u_V$, which is a particular case
of step 2 in the proof of theorem 4.3.4.

Let $u \in {\rm Hom}^\otimes(\eta,\omega)$ such that $f_u$ is a 
$*$-homomorphism. Then for every unitary $A$-comodule  $(V,\phi_V)$ we have
$u_{\overline
V} \circ \lambda_V = \overline u_V$, hence $\phi_V \circ (\overline u_V \otimes
u_V) = \phi_V \circ (u_{\overline V} \otimes u_V) \circ (\lambda_V \otimes
1_{\eta(V)}) = \phi_V \circ u_{\overline V \otimes V} \circ \widetilde
\eta_{\overline V,V} \circ (\lambda_V \otimes 1_{\eta(V)}) = u_{\mathbb C} \circ
\eta(\phi_V) \circ \widetilde \eta_{\overline V,V} \circ (\lambda_V \otimes
1_{\eta(V)}) = \widetilde \eta_1^{-1} \circ \eta(\phi_V) \circ \widetilde
\eta_{\overline V,V} \circ (\lambda_V \otimes 1_{\eta(V)}) = \phi'_V$, hence $u_V$ is unitary and
$u \in U^\otimes (\eta,\omega)$. \square

\subsection{Hopf $*$-co-Morita equivalence.}

We rewrite subsection 2.6 in the context of unitarizable Hopf $*$-algebras.

\begin{defi}
Let $A$ and $B$ be two unitarizable Hopf $*$-algebras.
Then $A$ and $B$ are said to be {\rm Hopf $*$-co-Morita equivalent} if there
is a monoidal $*$-functor $F : {\cal U}_f(A) \longrightarrow {\cal U}_f(B)$
with unitary monoidal constraints, which is an equivalence of categories.
\end{defi}

 Let $A$ be a unitarizable Hopf $*$-algebra and let
$\eta : {\cal U}_f(A) \longrightarrow {\cal H}$ be $*$-fibre functor. Then
the Hopf algebra ${\rm End}^\vee(\eta)$ can be endowed with a
Hopf $*$-algebra structure and a positive and faithful Haar measure (use the proof of proposition 4.3.3 or \cite{[W2]}). The unitarizable
Hopf $*$-algebras $A$ and
${\rm End}^\vee(\eta)$ are Hopf $*$-co-Morita equivalent.
Moreover $A$ and ${\rm End}^\vee(\eta)$ are isomorphic if and only if
there is a monoidal $*$-equivalence  (with unitary monoidal functor
constraints) $F :
{\cal U}_f(A) \rightarrow {\cal U}_f(A)$  such that the $*$-fibre
functors $\omega$ and $\eta\circ F$ are isomorphic ($\omega$ is the forgetful functor).

\section{$C^*$-norms on a Galois extension}

 Let $A$ be a unitarizable Hopf $*$-algebra and let $Z$ be an $A$-$*$-Galois extension endowed with a positive Haar measure. Theorem 5.2.1 ensures the existence of a $C^*$-norm on $Z$. This theorem is a generalization of theorem 4.4
in \cite{[DK]}. It allows us to find another proof of the Doplicher-Roberts'
unicity theorem for symmetric $*$-fibre functors. The notion of Galois extension for a compact quantum group is introduced. 

\subsection{Structure of Galois extensions}

\begin{theo}
 Let $A$ be a unitarizable Hopf $*$-algebra and let $Z$ be an $A$-$*$-Galois extension endowed with a positive Haar measure.
Let $\widehat A$ be the set
of isomorphism classes of irreducible comodules. For
$\lambda \in \widehat A$, let $n_\lambda$ be the dimension of a representative of
$\lambda$. Then there is a decomposition $Z = \bigoplus_{\lambda \in \hat
A} Z(\lambda)$, where $Z(\lambda)$ is a sub-$A$-comodule of dimension
$n_\lambda \times p_\lambda$ and  $p_\lambda$ is a positive integer.
For all $\lambda \in \widehat A$, there is a basis $(z_{ij})_{1 \leq i \leq
n_\lambda, 1 \leq j \leq p_\lambda}$ of $Z(\lambda)$ such that
$\sum^{n_\lambda}_{i=1} z^{\lambda^*}_{ij} z^\lambda_{ik} = \delta_{jk}$ and
$\sum^{p_{\lambda}}_{j=1} z^\lambda_{ij} z^{\lambda^*}_{kj} = \delta_{ik}$.
Moreover there is a basis $(a^{\lambda}_{ij})^{\lambda \in \widehat A}_{1 \leq i,j \leq n_\lambda}$ of $A$ such that each matrix $(a^{\lambda}_{ij})$ is 
unitary and $\alpha_Z(z^{\lambda}_{ij}) = \sum^{n_\lambda}_{k=1}  
a^{\lambda}_{ik} \otimes z^{\lambda}_{kj}$.
\end{theo}

\noindent
{\bf Proof}. By theorem 4.3.4 there is
 a $*$-fibre functor $\eta$ on
${\cal U}_f(A)$ such that ${\rm Hom}^\vee(\eta,\omega) \cong Z$. Thus by lemma 4.3.2 there is an isomorphism 
 $f :
\bigoplus_{\lambda \in \hat A} {\rm Hom}(\omega(V_\lambda), \eta(V_\lambda))
\longrightarrow Z$.
Let $Z(\lambda) = f({\rm Hom}(\omega(V_\lambda), \eta(V_\lambda))$.
The dimension of $Z(\lambda)$ is obviously equal to $n_\lambda \times p_\lambda$ where
$p_\lambda = \dim \eta(V_\lambda)$. Let us choose orthonormal bases
of $V_\lambda$ and $\eta(V_\lambda)$. The first equality was shown
in step 3 of the proof of theorem 4.3.4 and the second one can be proved in the same way. The proof of the last assertion is straightforward. \square

\begin{rema}
{\rm Let $A$ be Hopf algebra such that ${\rm Co}_f(A)$ is
semisimple and let $Z$ be an $A$-Galois extension. It is not
hard to show that there is a decomposition 
$Z = \bigoplus_{\lambda \in \widehat A} Z(\lambda)$
without using the  ${\rm Hom}^\vee$. It is even a simplified way to prove
Ulbrich's theorem in that case, avoiding  minimal models from \cite{[JS]}.
Then theorem 5.1.1 can be proved only using
the proof of proposition 4.3.1.}
\end{rema}

\subsection{Construction of a $C^*$-norm}

 Let $A$ be a unitarizable Hopf $*$-algebra and let $Z$ be an $A$-$*$-Galois extension endowed with a positive Haar measure.
Let $p$ be a $C^*$-semi-norm on $Z$. We have $p(z^\lambda_{ij}) \leq 1$ for all
the $z^\lambda_{ij}$'s of theorem 5.1.1. (see the first formula in this
theorem). Hence the upperbound of
$C^*$-semi-norms exists on $Z$.

Let us consider the regular representation $ L : Z \longrightarrow {\rm
End}(Z)$ defined by $L(x)(y) = L(xy)$. Let us note that $Z$ is prehilbert space, with $\langle x,y \rangle = \mu(x^*y)$ (the Haar measure $\mu$ is faithful by corollary 4.3.5). We deduce from the first formula of theorem 5.1.1 that $L(x)$ is continuous for all $x \in
Z$, and thus can be extended to a bounded operator on  $H$, the 
Hilbert space completion of $Z$. In this way we have a faithful $*$-representation $L' : Z
\longrightarrow {\cal B}(H)$ (where ${\cal B}(H)$ is the algebra of bounded linear operators
on $H$). We have proved a generalization of theorem 4.4 in
\cite{[DK]}:

\begin{theo}
 Let $A$ be a unitarizable Hopf $*$-algebra and let $Z$ be an $A$-$*$-Galois extension endowed with a positive Haar measure.
Then there is a $C^*$-norm on $Z$. Furthermore the upperbound of $C^*$-semi-norms exists on $Z$ (and hence is a $C^*$-norm).
\end{theo}

As an application of this result, we give a proof
 Doplicher-Roberts'
unicity theorem for symmetric $*$-fibre functors on the category
of unitary representations of a compact group.
Let $G$ be a compact group and let $R(G)$ be the unitarizable Hopf $*$-algebra 
of representative functions on $G$. Let
${\cal U}_f(G)$ be the monoidal $*$-category of finite-dimensional unitary representations of $G$. Obviously the categories ${\cal U}_f(G)$ and
${\cal U}_f(R(G))$ are identical. Moreover ${\cal U}_f(G)$ is endowed with a unitary symmetry (see \cite{[DR]}, 1.7) since $R(G)$ is commutative and 
thus ${\cal U}_f(R(G))$ is a symmetric monoidal $*$-category.
 A $*$-fibre functor on ${\cal U}_f(G)$ is said to be symmetric 
if it is symmetric as a monoidal functor.
We get (\cite{[DR]}, 6.9):

\begin{theo}
 Let $G$ be a compact group. Any symmetric $*$-fibre
functor on ${\cal U}_f(G)$ is (unitarily) isomorphic to the forgetful functor.
\end{theo}

\noindent
{\bf Proof}. Let $\eta$ be a symmetric $*$-fibre functor
and let $\omega$ be the forgetful functor. We must show that 
$U^\otimes(\eta,\omega) \not = \emptyset$. By theorem 4.4.1, 
$U^\otimes(\eta, \omega)$ and ${\rm Hom}_{*-{\rm alg}}({\rm
Hom}^\vee(\eta,\omega), \mathbb C)$ can be identified 
and by theorem 5.2.1 there is a $C^*$-norm on ${\rm Hom}^\vee(\eta,\omega)$.
But ${\rm Hom}^\vee(\eta,\omega)$ is commutative since
$\eta$ is symmetric. The Gelfand-Naimark theorem then ensures that
${\rm
Hom}_{*-{\rm alg}}({\rm Hom}^\vee(\eta,\omega),\mathbb C) \not =
\emptyset$. \square

\subsection{Galois extension for a compact quantum group}

We give a unusual definition of compact quantum group. 
However, it is equivalent to classical ones 
(\cite{[W3],[W1],[DK]}, see \cite{[MV]} for an expository text),
 thanks to the works of  Woronowicz and Van Daele (see below).

\begin{defi}
A {\rm compact quantum group} is a pair  $(A, ||.||)$ where $A$ is a Hopf $*$-algebra and $|| . ||$ is a $C^*$-norm on
$A$ such that the comultiplication is continuous. 
The Hopf $*$-algebra $A$ is called {\rm 
the algebra of representative function} on the compact quantum group. 
\end{defi}

Let  $(A, ||.||)$ be a compact quantum group. Let $\overline A$ be the $C^*$-algebra completion of $A$ and extend $\Delta$ to a $*$-homomorphism 
 $\overline \Delta : \overline A \longrightarrow \overline A
\bar \otimes \overline A$. 
The conditions of theorem 
2.4 in \cite{[V]} are fulfilled since $A$ is a Hopf algebra. Hence 
there is a state $\mu$ on $\bar A$ such that
$\mu * \varphi = \varphi * \mu = \mu$ where $\mu * \varphi = (\mu
\bar \otimes
\varphi) \circ \overline \Delta$. Then $\mu$ is a
Haar measure on $A$ and by \cite{[W1]}, 4.22, $\mu$ is faithful on
 $A$ : $\mu(a^*a)>0$
when $a\not = 0$ (see also \cite{[V2]}). Hence $A$ is
a unitarizable Hopf $*$-algebra.

Conversely, any unitarizable Hopf $*$-algebra is the algebra of representative 
functions on 
 a compact quantum group: see \cite{[DK]} or 5.2.1.

Our definition is a mix between the topological definition of Woronowicz
\cite{[W3]} and the algebraic one of Dijkhuizen an Koornwinder \cite{[DK]}.
It keeps some topological flavour but does not hide the group structure.

\begin{defi}
Let  $(A, ||.||)$ be a compact quantum group. 
An  {\rm $(A, ||.||)$-Galois extension} is a pair $(Z,||.||)$ where 
$Z$ is an $A$-$*$-Galois extension and $||.||$ is a $C^*$-norm on $Z$
such that the coaction $\alpha_Z$ is continuous.
\end{defi}

\begin{ex}
{\rm 1) Let $(A,||.||)$ be a compact quantum group. Then
$(A,||.||)$ is an $(A, ||.||)$-Galois extension.

\noindent
2) Let $(A,||.||)$ be a compact quantum group and let $Z$ be an
 $A$-$*$-Galois extension endowed with a positive Haar measure.
 Let $||.||_\infty$ be the upperbound of $C^*$-semi-norms on $Z$.
Then $(Z, ||.||_\infty)$ is an $(A,||.||)$-Galois extension.}
\end{ex}

Actions of compact quantum groups on quantum spaces were considered
by Podles in \cite{[Po]}. Theorem 1.5 in \cite{[Po]} shows that our definition
of actions (with the existence of a dense comodule-algebra)
is the same as the one given in \cite{[Po]}.
Measure theory can be deduced from topology:

\begin{theo}
 Let $(A,||.||)$ be a compact quantum group and let 
$(Z,||.||)$ be an $(A,||.||)$-Galois extension. Let  $\overline A$
(resp.  $\overline Z$) be the $C^*$-algebra completion of $A$ (resp. $Z$).
 There is a state $\mu$ on  $\overline Z$
such that for every state $\varphi$ on $\overline A$, then
$\varphi*\mu=\mu$ (where 
$\varphi*\mu= (\varphi \bar \otimes \mu) \circ \overline
\alpha_Z$ and $\alpha_Z$ has be extended  in $\overline
\alpha_Z : \overline Z \longrightarrow \overline A \bar \otimes \overline
Z$). Moreover if $z\in Z$ is a positive (as an element 
the $C^*$-algebra $\overline Z$) non-zero element, then $\mu(z)>0$.
In particular $\mu$ is a positive Haar measure on $Z$.
\end{theo}

\noindent
{\bf Proof}. Let us first remark that there is a faithful Haar measure on
$Z$ by proposition 4.2.5 and corollary 4.3.5. But the assertion in the theorem is more precise.

Let $J$ be the Haar measure on $\overline A$ given by
theorem 2.4 in \cite{[V]} and let $\psi$ be a state on $\overline Z$. 
Then $\mu =
J*\psi$ is a positive Haar measure on $Z$ by proposition 4.2.5, and does not 
depend on the choice of the state $\psi$ by proposition 4.2.2. Clearly
$\mu$ satisfies the first claim in the theorem.

We follow the ideas of the proof of (4.22) in \cite{[W1]}.
Let $z\in Z$ be a positive element of $\overline Z$ such that
$\mu(z)=0$. Then for every state $\rho$ on $\overline Z$, 
we have that $z* \rho =
1_A\otimes \rho(\alpha_Z(z))\in A$ is a positive element of 
$\overline A$. But $J(z* \rho)=J*\rho(z)= \mu(z)=0$ ($J*\rho=\mu$), 
therefore $z*\rho = 0$ by \cite{[W1]},
4.22. Every continuous linear functional on $\overline Z$ is a linear
combination of states and hence $z*\rho=0$ for
every continuous linear functional $\overline Z$.

Let us use theorem 5.1.1 : $Z=\oplus_{\lambda \in \hat A}
Z(\lambda)$  with  basis $(z_{ij}^\lambda)$ and
$A=\oplus_{\lambda \in \hat A} A(\lambda)$ with basis $(a_{ij}^\lambda)$, 
such that $\alpha_Z(z_{ij}^\lambda) = \sum_k a_{ik}^\lambda \otimes
z_{kj}^\lambda$.
Let $F$ be a finite subset of $\widehat A$ and let $\xi_{ij}^\lambda$ 
be complex numbers such that
$z=\sum_{\lambda \in F} \sum_{i,j} \xi_{ij}^\lambda
z_{ij}^\lambda$. Then $z*\rho = \sum_{\lambda \in F} \sum_{i,j,k}
\xi_{ij}^\lambda \rho(z_{kj}^\lambda) a_{ik}^\lambda = 0$
and hence $\sum_j \xi_{ij}^\lambda \rho(z_{kj}^\lambda) = 0$
for all $\lambda \in F$ and all $i$ and $k$
(the $a_{ij}^\lambda$'s  are linearly independent).
This is true for every continuous linear functional
 $\rho$ on $\overline Z$. But the $(z_{ij}^\lambda)$'s
are linearly independent, and therefore the Hahn-Banach theorem ensures that
$\xi_{ij}^\lambda=0$ for all $i$ and $j$, which means that $z=0$.
\square

\section{Universal Galois extensions}

In this section we introduce universal Galois extension. They are useful
to ensure the existence of a positive Haar measure on a Galois extension.
As an application, the example found in 2.7 is adapted to the compact quantum 
group setting and gives rise to a Galois extension for the compact quantum group $U_q(2)$.

\subsection{The $C^*$-algebra generated by an unitary matrix}

We first define an important class of $C^*$-algebras. 

Let us fix some notations.
Let $A$ be a  $*$-algebra, and let 
$a = (a_{ij})_{1 \leq i \leq n, 1 \leq j \leq
p} \in M_{n,p}(A)$. We define $\overline a = (a^*_{ij})_{1 \leq i \leq n, 1
\leq j \leq p}$ and $a^* = (a^*_{ji})_{1 \leq j \leq p, 1 \leq i \leq n} \in
M_{p,n}(A)$. We will say that $a$ is unitary if $a^*a = I_{M_{p(A)}}$
 and $aa^* = I_{M_n(A)}$ and we will write $a^*a = aa^* = 1$.

\begin{defi}
The $*$-algebra $O^o_{n,p}$ is the universal $*$-algebra 
generated by a unitary  $(n,p)$-matrix:
$
O^o_{n,p} \ \hbox{\rm is the quotient of} \ \mathbb C\{a_{ij},a^*_{ij}\}_{1 \leq i
\leq n, 1 \leq j \leq p}
$
by the two-sided  $*$-ideal generated by the  $2np$ relations $aa^* = a^*a = 1$.
\end{defi}

It is easily seen that there is a non-trivial $*$-representation of $O^o_{n,p}$
on a separable Hilbert space, and hence a $C^*$-semi-norm on  $O^o_{n,p}$. 

The upperbound of $C^*$-semi-norms on $O^o_{n,p}$ exists since the matrix $a$
is unitary . For all $x \in O^o_{n,p}$, let
$||x||_\infty = \sup ||\pi(x)||$ where $\pi$ runs over all Hilbert space 
$*$-representations of $O^o_{n,p}$. Let $N = \{x \in
O^o_{x,p}, ||x||_\infty = 0\}$. Then $N$ is a two-sided $*$-ideal of 
$O^o_{n,p}$ and  $|| . ||_\infty$ produces a $C^*$-norm on $O^o_{n,p}
/ N$. 

\begin{defi}
The  $C^*$-algebra $O_{n,p}$ is the $C^*$-algebra completion of $O^o_{n,p}/N$ with respect to $||.||_\infty$.
\end{defi}

Universal property of $O_{n,p}$: there is an obvious $*$-homomorphism  $i :
O^o_{n,p} \longrightarrow O_{n,p}$ and if  $\pi : O^o_{n,p} \longrightarrow B$ is a $*$-homomorphism into a $C^*$-algebra $B$, there is a unique 
$*$-homomorphism $\tilde \pi : O_{n,p} \longrightarrow B$ such that
$\tilde \pi \circ
i= \pi$.

\smallskip
There is a $*$-algebra isomorphism $\varphi : O^o_{n,p} \longrightarrow
O^o_{p,n}$ given by $\varphi(a_{ij}) = a^*_{ji}$, and thus by the 
universal property the $C^*$-algebras $O_{n,p}$ and $O_{p,n}$ are isomorphic.

When $p = 1$ the $C^*$-algebra  $O_{n,1}$ is isomorphic to the Cuntz algebra
$O_n$ (\cite{[C]}). If $n \not = 1$, then $O_n$ is a simple $C^*$-algebra
 (\cite{[C]}). The following question arises naturally: 
if  $n \not = p$, is $O_{n,p}$ 
a simple $C^*$-algebra?

\subsection{Construction of universal Galois extensions}

\begin{defi}
A Hopf $*$-algebra $A$ is said to be a  
{\rm matrix Hopf $*$-algebra} if it is generated (as a $*$-algebra) by
elements $(a_{ij})_{1 \leq i,j \leq n}$ such that $\Delta(a_{ij}) = \sum a_{ij}
\otimes a_{jk}$ and  $\varepsilon (a_{ij}) = \delta_{ij}$.

A {\rm compact matrix quantum group} is a compact quantum group $(A,||.||)$ such that $A$ is a  matrix Hopf $*$-algebra
\end{defi}

\begin{theo}
 Let $A$ be a unitarizable matrix Hopf $*$-algebra and let
 $Z$ be an $A$-$*$-Galois extension endowed with a positive Haar measure.
Then there is a surjective (and non injective) $*$-homomorphism
$\varphi : O^o_{n,p}
\longrightarrow Z$ for some integers $n$ and $p$.
Moreover there is a surjective $*$-homomorphism
 $\tilde \varphi : O_{n,p} \longrightarrow \overline Z$
where $\overline Z$ denotes the $C^*$-algebra completion of
$Z$ with respect to an arbitrary  $C^*$-norm.
\end{theo}

\noindent
{\bf Proof}. Let $\eta$ be a
 $*$-fibre functor  on
${\cal U}_f(A)$ such that ${\rm Hom}^\vee(\eta,\omega) \cong Z$.
Let $(a_{ij})_{1 \leq i, j \leq n}$ be generators of $A$ and let 
$V$ be the obvious $n$-dimensional comodule associated to the $(a_{ij})$'s. 
There are invariant scalar products on $V$ and  $\overline V$, and
for every object 
$W$ of ${\cal U}_f(A)$, there is an isometry from $W$ into a direct sum
of tensor powers of $V$ and $\overline V$. 
This fact implies that ${\rm Hom}^\vee(\eta,\omega)$ and therefore  
$Z$ is generated 
as a $*$-algebra by elements $(z_{ij})_{1\leq i \leq n, 1 \leq j \leq p}$
where $p =\dim \eta(V)$. 
Moreover one can assume that these elements satisfy the relations 
of $O^o_{n,p}$ (choose orthonormal bases of $V$ and $\eta(V)$ and use step 
3 in the proof of theorem 4.3.4). We get a surjective
 $*$-homomorphism $\varphi : O^o_{n,p} \longrightarrow Z$.

The comodule $\overline V$ is unitarizable. Let us choose orthonormal bases
of $\overline V$ and $\eta(\overline V)$ : it is easily seen that there are 
matrices $G \in GL_p(\mathbb C)$ and $F \in GL_n(\mathbb C)$ 
such that the matrix
$F\overline z G^{-1}$ is unitary 
(the matrix $F$ is such that 
$F\overline aF^{-1}$ is  unitary).
This means that $\varphi$ is not injective. 

Let $\tilde \varphi : O_{n,p}
\longrightarrow \overline Z$ be the $C^*$-homomorphism produced by
the universal property of $O_{n,p}$. Then $\tilde \varphi$ is surjective 
since  $\tilde
\varphi(O_{n,p})$ is  a sub-$C^*$-algebra of $\overline Z$ containing
a dense sub-algebra. \square

\begin{rema}
{\rm We could say that $\tilde \varphi$ is not injective if we knew
that  the ideal $N$ of 6.1 is trivial.}
\end{rema}

The proof of this result leads us to the following definition:

\begin{defi}
Let $G \in GL_p(\mathbb C)$ and let $F \in GL_n(\mathbb C)$.
The $*$-algebra
$A^o_u(F,G)$ is the quotient of the free $*$-algebra $\mathbb
C\{z_{ij},z^*_{ij}\}_{1 \leq i
\leq n, 1
\leq j
\leq p}$ by the two-sided $*$-ideal generated by the relations $zz^* = z^*z = 1$ and $F
\overline z G^{-1}$ unitary.
\end{defi}

When $n = p$ and $F = G$, the above defined algebras are the universal
quantum groups of Van Daele and Wang (\cite{[VW]}): $A^o_u(F) = A^o_u(F,F)$ 
is an unitarizable Hopf $*$-algebra 
(by its definition), whose representation theory  is studied in
\cite{[Ba2]}. The following result is contained in the proof of theorem 6.2.2.

\begin{coro}

 Let $A$ be a unitarizable matrix Hopf $*$-algebra and let
 $Z$ be an $A$-$*$-Galois extension endowed with a positive Haar measure.
Then there is a surjective $*$-homomorphism
$\varphi : A^o_u(F,G) \longrightarrow Z$ for some matrices
$G \in GL_p(\mathbb C)$ and $F \in GL_n(\mathbb C)$.
\end{coro}

The next result is an useful criterion to ensure the existence of a 
positive Haar measure.

\begin{prop}

Let $A$ be a unitarizable matrix Hopf $*$-algebra 
with  $\varphi : A^o_u(F) \longrightarrow A$ a surjective $*$-homomorphism
for some
$F \in GL_n(\mathbb C)$.
Let $Z$ be an
$A$-$*$-Galois extension. 
Assume that there is a surjective $*$-homomorphism
$\psi : A^o_u (F,G) \longrightarrow Z$ for some 
$G \in GL_p(\mathbb C)$. Then there is a positive Haar measure on
$Z$.
\end{prop}

\noindent
{\bf Proof}. Let $(a_{ij})_{1 \leq i,j \leq n}$ be the generators of $A$ 
provided by $\varphi$: we have $a^* a
= aa^* = 1$ and the matrix $F\overline a F^{-1}$ is unitary.
Let $V$ be the obvious $n$-dimensional unitary $A$-comodule associated to
the $a_{ij}$'s with orthonormal basis $v_1,\ldots,v_n$. 
Let $(z_{ij})_{1 \leq i \leq n, 1 \leq j
\leq p}$ be the generators of $Z$ provided by $\psi$: we have
$zz^* = z^*z = 1$ and $F\overline z G^{-1}$ is unitary. 

By theorem 4.3.4 it is sufficient to show that
the fibre functor $\eta_Z : {\cal U}_f(A)
\longrightarrow {\rm Vect}_f(\mathbb C)$ of 2.2
 factorizes through a $*$-fibre functor
$\eta^\odot_Z : {\cal U}_f(A) \longrightarrow {\cal H}$. 

We  first show that the maps defined in step 1 of the proof of proposition 4.3.1 are scalar products on $V\wedge Z$ and $\overline V \wedge Z$.
Let $\sigma_j = \Sigma_i v_i \otimes z_{ij}$. It is obvious that 
$(\sigma_j)$ is an orthonormal basis for $V\wedge Z$.
Now let $t = (t_{ij}) = F \overline z G^{-1}$ and  let
$\overline w_1,
\ldots, \overline w_n$ be an orthonormal basis of $\overline V$.
It is easy to check that the elements $\tau_j = \sum_i \overline \omega_i \otimes t_{ij}$ belong to
$\overline V \wedge Z$ and that $(\tau_j)$ is an orthonormal basis
for the sesquilinear map defined in step 1 of the proof of proposition 4.3.1.

Let $\mathcal C$ be the full subcategory of  ${\cal U}_f(A)$ whose objects
are tensor powers of  $V$ and $\overline V$. It is then easily seen that
the sesquilinear maps defined in step 1 of proposition 4.3.1 are scalar
products for all objects of $\mathcal C$ and that the monoidal constraints
are unitary maps. Thus we get a monoidal $*$-functor 
$\eta^\odot_Z : \mathcal C \longrightarrow {\cal H}$.
Every object of  ${\cal U}_f(A)$ can be isometrically embedded into
a direct sum of objects of $\mathcal C$. Hence $\eta^\odot_Z$
can be extended into a $*$-fibre functor 
$\eta^\odot_Z : {\cal U}_f(A) \longrightarrow {\cal H}$ (see lemma 6.5
in \cite{[DR]}, the whole $*$-structure of ${\cal U}_f(A)$ is determined by
its isometries). \square 

\bigskip

 Let $G \in GL_p(\mathbb C)$ and let $F \in GL_n(\mathbb C)$. Assume that
$A^o_u (F,G)$ is a non-zero algebra. Then $A^o_u (F,G)$ is an
$A^o_u (F)$-$*$-Galois extension. This can be proved in the same way
we have shown that $GL_q(2,-2)$ is a $GL_q(2)$-Galois extension.
By proposition 6.2.6 there is a positive Haar measure on
$A^o_u (F,G)$ and therefore there is a $C^*$-norm on $A^o_u (F,G)$ by
theorem 5.2.1. Let us denote by $A_u(F,G)$ the $C^*$-algebra completion
of $A^o_u (F,G)$ with respect to the maximal $C^*$-norm.

\begin{coro}
Let $(A,||.||)$ be a compact matrix quantum group and let $(Z,||.||)$ be 
an $(A,||.||)$-Galois extension. Let $\overline Z$ be the $C^*$-algebra
completion of $Z$. Then there is a surjective $*$-homomorphism
$\tilde \psi : A_u (F,G) \longrightarrow \overline Z$ for some
matrices  $G \in GL_p(\mathbb C)$ and $F \in GL_n(\mathbb C)$.
\end{coro}

Corollaries 6.2.5 and 6.2.7 justify the term universal Galois extension
for our algebras $A^o_u(F,G)$ and $A_u(F,G)$.

\medskip

We end this subsection by a discussion on the preservation of dimensions by a
$*$-fibre functor.
Let $A$ be an unitarizable Hopf $*$-algebra and let $V$ 
be a $n$-dimensional unitary $A$-comodule.
Let $\eta$ be a $*$-fibre functor on  ${\cal U}_f(A)$
such that $\dim \eta(V) = p$ where $p \not = n$. Let $\langle V \rangle$ be
the full sub-category of ${\cal U}_f(A)$ whose objects are sub-objects
of direct sums of tensor powers of $V$ and $\overline V$. Then
$\langle V \rangle$ is a monoidal $*$-category with conjugation, 
and there is a unitarizable matrix Hopf $*$-algebra $B$ such that
 ${\cal U}_f(B)
{\widetilde {\longrightarrow}} \langle V \rangle$ 
(\cite{[W2]}). 
Hence $\eta$ induces a $*$-fibre functor on 
${\cal U}_f(B)$ and we get a 
$B$-$*$-Galois extension $Z$ endowed with a positive Haar measure.
By theorem 6.2.2 there is a surjective (and non-injective) $*$-homomorphism
$O^o_{n,p} \longrightarrow Z$. If we knew that
every Hilbert space $*$-representation of $O^o_{n,p}$ is faithful (this is
true in the Cuntz algebra case $p=1$), we could conclude that $n=p$ and
that a $*$-fibre functor preserves the dimensions of the underlying vector spaces.

\subsection{A Galois extension for $U_q(2)$}

We use the notation and results of 2.7. We now have $k=\mathbb C$ and 
$q\in {\mathbb R}-\{0\}$.

Let us denote by $U_q(2)$ the Hopf $*$-algebra whose underlying Hopf algebra
is $GL_q(2)$ and whose $*$-structure is given by
$x^*_{11} = x_{22}t$, $x^*_{12} =
-q^{-1}x_{21}t$, $x^*_{21} = -qx_{12}t$, $x^*_{22} = x_{11}t$ and $t^* =
z_{11}z_{22} - q^{-1}z_{12} z_{21}$. It is well known that $U_q(2)$
is unitarizable.

\begin{prop}
There is a $*$-algebra structure on $GL_q(2,-2)$ such that $z^*_{11} =
z_{22} \tau$, $z^*_{12} = q^{-1}z_{21} \tau$, $z^*_{21} = q\tau z_{12}$,
$z^*_{22} =
\tau z_{11}$ et $\tau^* = z_{11} z_{22} + q^{-1}z_{12} z_{21}$.
Let us denote by $U_q(2,-2)$ the associated $*$-algebra.
The coaction defined in lemma 2.7.2 turns  $U_q(2,-2)$ into a
$U_q(2)$-$*$-Galois extension. Furthermore the Haar measure on
 $U_q(2,-2)$ is positive.
\end{prop}

\noindent
{\bf Proof}. We leave all the computations to the reader. The only statement
which does not follow from routine computations is the positivity of the Haar 
measure. Let
 $F_q\in GL_3(\mathbb C)$ be the matrix 
$F_q =  \left(\begin{array}{ccc} 0 & 1 & 0 \\
-q & 0 & 0 \\
0 & 0 & 1 \end{array} \right)$.
There is a surjective $*$-homomorphism
$\varphi :
A^o_u(F_q,F_{q}) \longrightarrow U_q(2)$ defined by 
$\varphi(a_{ij}) = x_{ij}$ if  
$i$ et $j \leq 2, \varphi(a_{33}) = t $ and  $\varphi(a_{ij}) = 0$ if
$i>2$ or $j>2$
(in particular $U_q(2)$ is unitarizable). In the same way there is 
a surjective $*$-homomorphism $\psi : A^o_u(F_q,F_{-q}) \longrightarrow 
U_q(2,-2)$ defined by the same formulas. By proposition 6.2.6
there is a positive Haar measure on $U_q(2,-2)$. \square

\bigskip
Let  $||.||_\infty$ be the upperbound of $C^*$-semi-norms on 
$U_q(2)$ and $U_q(2,-2)$ respectively. Proposition 6.3.1 is translated in :

\begin{coro}
$(U_q(2,-2), ||.||_\infty)$ is a $(U_q(2),||.||_\infty)$-Galois extension.
\end{coro}

\bigskip

D\'epartement des Sciences Math\'ematiques,
case 051

Universit\'e Montpellier II

Place Eug\`ene Bataillon, 34095 Montpellier Cedex 5

{\tt e-mail : bichon\char64math.univ-montp2.fr}

\end{document}